\documentclass[11pt]{article}
\usepackage{amsfonts, amsmath, amssymb, amsthm, float, subcaption, setspace, enumitem}

\usepackage[T1]{fontenc}
\usepackage[utf8]{inputenc}

\usepackage[square,numbers]{natbib} 
\usepackage{graphicx}

\usepackage[a4paper, margin=2.4cm]{geometry}
\usepackage{xcolor} 

\usepackage{tikz} 
\usetikzlibrary{shapes.geometric}

\usepackage{hyperref,cleveref} 

\definecolor{theme}{RGB}{170, 0, 0}

\hypersetup{
    colorlinks=true,
    linkcolor=black, 
    filecolor=black,
    citecolor=black, 
    urlcolor=black,
    pdftitle={}, 
    pdfpagemode=FullScreen,
    }

\DeclareMathOperator{\tw}{tw}

\DeclareMathOperator{\bigo}{\mathcal{O}}
\renewcommand{\leq}{\leqslant}
\renewcommand{\geq}{\geqslant}

\theoremstyle{plain} 
\newtheorem{theorem}{Theorem}
\newtheorem{lemma}[theorem]{Lemma}

\newtheorem{observation}[theorem]{Observation}
\newtheorem{claim}[theorem]{Claim}

\newtheorem{problem}[theorem]{Problem}

\theoremstyle{definition} 
\newtheorem{definition}[theorem]{Definition}

\newcommand{\define}[1]{{\color{theme}\textit{#1}}} 

\begin{document}
\renewcommand{\thefootnote}{\fnsymbol{footnote}} 

\title{\bf On Universal Graphs for Trees and Tree-Like Graphs} 

\author{
Neel Kaul\footnotemark[1] \and
Jaehoon Kim\footnotemark[2] \and  
Minseo Kim\footnotemark[2] \and 
David R.~Wood\footnotemark[4]
}

\footnotetext[1]{School of Mathematics, Monash University, Melbourne, Australia (\texttt{neel.kaul@monash.edu}).} 

\footnotetext[2]{Department of Mathematical Sciences, KAIST, South Korea (\texttt{\{jaehoon.kim, minseo.kim00\}@kaist.ac.kr}).     Jaehoon Kim supported by the National Research Foundation of Korea (NRF) grant funded by the Korea government(MSIT) No. RS-2023-00210430.}

\footnotetext[4]{School of Mathematics, Monash University, Melbourne, Australia (\texttt{david.wood@monash.edu}). Research supported by the Australian Research Council and by NSERC. } 

\maketitle

\begin{abstract}
    Chung and Graham~[\textit{J.~London Math.\ Soc.\ }1983] claimed to prove that there exists an $n$-vertex graph $G$ with \( \frac{5}{2}n \log_2 n + \bigo(n)\) edges that contains every $n$-vertex tree as a subgraph. Frati, Hoffmann and T\'oth~[\textit{Combin.\ Probab.\ Comput.} 2023] discovered an error in the proof. By adding more edges to $G$ the error can be corrected, bringing the number of edges in $G$ to \(\frac{7}{2}n \log_2 n + \bigo(n). \) 

    We make the first improvement to Chung and Graham's bound in over four decades by showing that there exists an $n$-vertex graph with \(\frac{14}{5}n \log_2 n + \bigo(n) \) edges that contains every $n$-vertex tree as a subgraph.

    Furthermore, we generalise this bound for treewidth-$k$ graphs by showing that there exists a graph with $\bigo(kn\log(n/k+1))$ edges that contains every $n$-vertex treewidth-$k$ graph as a subgraph. This is best possible in the sense that $\Omega(kn\log(n/k+1))$ edges are required.
\end{abstract}

\renewcommand{\thefootnote}{\arabic{footnote}} 
\section{Introduction}

A graph $G$ is \define{$\mathcal{F}$-universal} for a graph family $\mathcal{F}$ if $G$ contains\footnote{A graph $G$ \define{contains} a graph $H$ if a subgraph of $G$ is isomorphic to $H$.} every graph in $\mathcal{F}$. The minimum number of edges in an $\mathcal{F}$-universal graph has been extensively studied for various graph families $\mathcal{F}$.
For example, families of trees \cite{bergold2024, grahamfirst, grahamthird, grahamsecond, pipenger, montgomery}, families of separable bounded degree graphs \cite{chung}, families of bounded degree graphs \cite{alonsecond, szemeredi, alonthird, capalbosecond}, families of $m$-edge graphs \cite{alon, babai, bucic}, and families of planar graphs with certain restrictions \cite{bhatt, capalbo, joret}. 

Let \define{$\mathcal{T}_n$} be the family of $n$-vertex trees. \citet{grahamfirst, grahamthird} considered the natural problem of determining the minimum number of edges in a $\mathcal{T}_n$-universal graph, denoted by \define{$s^{*}(n)$}. They further defined a more restrictive parameter, \define{$s(n)$}, as the minimum number of edges in an $n$-vertex $\mathcal{T}_n$-universal graph. Naturally, $s^{*}(n)\leq s(n)$. \citet{grahamfirst} proved the lower bound $s^*(n)\geq\frac{1}{2}n\ln{n}$. Furthermore, they claim to prove the upper bound $s(n)\leq\frac{5}{2}n\log_2{n}+\bigo(n)$ in \cite{grahamthird}, implying that $s^{*}(n)$ and $s(n)$ are $\Theta(n\log n)$. 

\citet{frati} discovered an error in the proof of the upper bound in \cite{grahamthird} and provided a fix at the cost of a slightly weaker bound. The same error was later independently rediscovered by the authors of this paper.\footnote{This paper merges and extends the results of \cite{Kaul, KimKim}.} In \Cref{sec: 3}, we explain the error in Chung and Graham's original paper and describe how to patch the error to obtain the bound,\footnote{Curiously, \citet*{bhatt} and \citet*{chung} both cited the upper bound from \cite{grahamthird} as $\frac{7}{2}n\log_2 n + \bigo(n)$.}
\begin{align}\label{eqn: corrected bound}
s(n)\leq\tfrac{7}{2}n\log_2{n}+\bigo(n)\;.
\end{align}

Recently, \citet{gyori} improved the lower bound to $s^{*}(n)\geq n\ln{n}-\bigo(n)$. So the best previously known bounds for $s^{*}(n)$ and $s(n)$ are:
\begin{align}\label{intro: ineq 1}
    n\ln{n}-\bigo(n)\leq s^{*}(n)\leq s(n)\leq \tfrac{7}{2}n\log_2{n}+\bigo(n)\;.
\end{align}
Note that the multiplicative gap between the two bounds is $\frac{7}{2\ln2}+o(1) \approx5.04$. In \Cref{sec:4}, we prove a new upper bound.
\begin{theorem}\label{thm2}
    $s(n)\leq \frac{14}{5}n\log_2{n} + \bigo(n)$.
\end{theorem}
This bound yields a new multiplicative gap between the lower and upper bounds of $\tfrac{14}{5\ln 2} \approx 4.04$. This is the first improvement to the upper bound in \labelcref{eqn: corrected bound} in over four decades. It is an interesting open problem to determine the optimal constant in the $n\log n$ term; we know
\begin{align*}
     n\ln{n}-\bigo(n)\leq s^{*}(n)\leq s(n)\leq \tfrac{14}{5\ln 2}n\ln{n}+\bigo(n)\;.
\end{align*}
We generalise these results to tree-like graphs through the notion of treewidth, which is a graph parameter that measures how close a given graph is to being a tree, and is of particular importance in graph minor theory and in algorithmic graph theory. See \Cref{def: treewidth} for the definition. 

Let \define{$\mathcal{T}^{k}_{n}$} be the family of $n$-vertex treewidth-$k$ graphs. Let \define{$s_{k}^{*}(n)$} be the minimum number of edges in a $\mathcal{T}^{k}_{n}$-universal graph. The best previously known upper bound for $s^{*}_{k}(n)$ is $\bigo(k n\log^{2}n)$ as explained by \citet{JoretTalk}. We determine $s^{*}_{k}(n)$ up to a constant factor. 
\begin{theorem}\label{thm: treewidth}
    $s_{k}^{*}(n)=\Theta(k n\log (n/k+1))$ whenever $n>k$.
\end{theorem}
In particular, we show that if $n > k$, then there exists a $\mathcal{T}^{k}_{n}$-universal graph with $\bigo(k n\log (n/k+1))$ edges, $\bigo(n)$ vertices and treewidth $\bigo(k\log (n/k+1))$. The treewidth of the universal graph is optimal in the sense that every $\mathcal{T}^{k}_{n}$-universal graph with $\bigo(n)$ vertices has treewidth $\Omega(k \log (n/k+1))$, assuming $n>k$.

\section{Preliminaries}\label{ch2}

The number of vertices in a graph $G$ is denoted by \define{$|G|$}. The \define{neighbourhood} of a vertex $u$ in $G$ is denoted by \define{$N_{G}(u)$}. A graph $H$ \define{embeds} into a graph $G$ if there is an injective homomorphism $\phi:V(H)\rightarrow V(G)$ (called an \define{embedding}). 

Consider an $n$-vertex tree $T$ rooted at a vertex $v$.
For a vertex $u$ of $T$, consider the unique path $(u_0,u_1,\ldots,u_t)$ where $u_{0}=u$ and $u_{t}=v$. The \define{level} of $u$ in $T$, denoted by \define{$L_{T}(u)$}, is the number $t$. A \define{cousin} of $u$ is a vertex with the same level as $u$. The \define{level} of $T$, denoted \define{$L(T)$}, is the maximum level of a vertex of $T$. The \define{$i$-th ancestor} of $u$ is $u_{i}$, and $u$ is an \define{$i$-th descendant} of $u_{i}$. Note that $u$ is the $0$-th ancestor and $0$-th descendant of itself. For $i>t$, we define the $i$-th ancestor of $u$ to be the root of $T$. The \define{parent} of $u$, denoted by \define{$u^{*}$}, is the vertex $u_{1}$ (provided $u$ is not the root of $T$), and $u$ is a \define{child} of $u^{*}$. A \define{sibling} of $u$ is a vertex with the same parent as $u$. 
The \define{closest common ancestor} of two vertices $u$ and $w$ of $T$
is the unique common ancestor of $u$ and $w$ with maximum level in $T$. Let \define{$D_{T}[u]$} be the set of all vertices $w$ of $T$ such that $w$ is an $i$-th descendant of $u$ for some $i\geq 0$.
Note that $u\in D_{T}[u]$.
Let $\define{$\nu_T(u)$}:=|D_T[u]|$ and let $\define{$D_{T}(u)$}:=D_{T}[u]\setminus\{u\}$. We omit the subscript $T$ from the notation when it is clear from the context.

Further assume that for each vertex $u$ of $T$, its children are ordered from left to right. Then we obtain a linear order on $V(T)$ by DFS preorder traversal\footnote{Starting at the root of $T$, visit the left-most unvisited child of the current vertex $w$. If every child of $w$ has been visited, backtrack to the parent of $w$. This gives a walk in $T$ that starts and ends at the root; the vertices of $T$ are ordered by their first appearance.}, say $w_{1},w_{2},\dots,w_{n}$. A rooted tree is \define{DFS-ordered} if its vertices are ordered by DFS preorder traversal $w_{1},w_{2},\dots,w_{n}$ from left to right. A subset $S$ of $V(T)$ is \define{$T$-admissible} if $S=\{w_{1},\dots,w_{t}\}$ for some $t\leq n$. Note that $\varnothing$ and $\{w_{1}\}$ are $T$-admissible. A \define{left-sibling}/\define{left-cousin} of a vertex $u$ of $T$ is a sibling/cousin of $u$ that comes before $u$ in the DFS-order. A \define{right-sibling}/\define{right-cousin} of $u$ is a sibling/cousin of $u$ that comes after $u$ in the DFS-order. The \define{nearest left-cousin} of $u$, denoted by \define{$l(u)$}, is the right-most left-cousin of $u$ in the DFS-order (provided $u$ has a left-cousin). Note that if $u$ is not the left-most child of its parent, then $l(u)$ is a left-sibling of $u$, otherwise, $u$ is the left-most child of its parent and $l(u)$ is the right-most child of $l(u^{*})$. The \define{nearest right-cousin} of $u$, denoted by \define{$r(u)$}, is the left-most right-cousin of $u$ in the DFS-order (provided $u$ has a right-cousin).

For a DFS-ordered tree $T$, we define a directed (multi)graph on the vertex set $V(T)$ as follows. For integers $g$ and $h$, let $\define{$[g, h]$}:=\{g, \ldots, h\}$ when $g\leq h$, and let $\define{$[h]$}:=\{1,\ldots,h\}$.

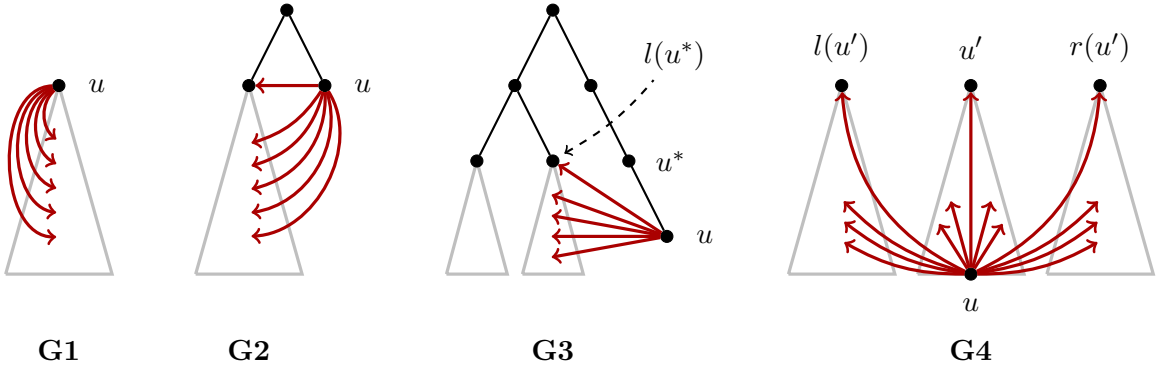
\begin{figure}[!ht]
    \centering
    \begin{tikzpicture}[/tikz/xscale=1, /tikz/yscale=1]
        \begin{scope}
            \draw[very thick, draw=lightgray] (-0.7,-0.5) -- (0.7,-0.5) -- (0,2) -- (-0.7,-0.5);
            \node (a) at (0.1,1.2) {};
            \node (b) at (0.1,0.9) {};
            \node (c) at (0.1,0.6) {};
            \node (d) at (0.1,0.3) {};
            \node (e) at (0.1,0) {};
        
            \node(index) at (0,-1.5) {\textbf{G1}};
        \end{scope}
        
        \begin{scope}
            \node[draw, shape = circle, fill = black, minimum size = 0.1cm, inner sep=1.6pt] (u) at (0,2) {};
            \node (text) at (0.5,2) {$u$};
        \end{scope}
        
        \begin{scope}[every edge/.style={draw=theme,very thick}]
            \path [->] (u) edge[bend right=50] (a);
            \path [->] (u) edge[bend right=60] (b);
            \path [->] (u) edge[bend right=70] (c);
            \path [->] (u) edge[bend right=80] (d);
            \path [->] (u) edge[bend right=90] (e);
        \end{scope}
    
        \begin{scope}[shift={(4-1.5,0)}]
            \draw[very thick, draw=lightgray] (-0.7,-0.5) -- (0.7,-0.5) -- (0,2) -- (-0.7,-0.5);
            \node (a) at (-0.1,1.2) {};
            \node (b) at (-0.1,0.9) {};
            \node (c) at (-0.1,0.6) {};
            \node (d) at (-0.1,0.3) {};
            \node (e) at (-0.1,0) {};
        
            \node[draw, shape = circle, fill = black, minimum size = 0.1cm, inner sep=1.6pt] (sibling) at (0,2) {};
            \node (text sibling) at (-0.7,2) {};
        
            \node[draw, shape = circle, fill = black, minimum size = 0.1cm, inner sep=1.6pt] (root) at (0.5,3) {};
        
            \node(index) at (0,-1.5) {\textbf{G2}};
        \end{scope}
        
        \begin{scope}[shift={(4-1.5,0)}]
            \node[draw, shape = circle, fill = black, minimum size = 0.1cm, inner sep=1.6pt] (u) at (1,2) {};
            \node (text) at (1.5,2) {$u$};
        \end{scope}
        
        \begin{scope}[every edge/.style={draw, thick}]
            \path [-] (root) edge (u);
            \path [-] (root) edge (sibling);
        \end{scope}
        
        \begin{scope}[every edge/.style={draw=theme,very thick}]
            \path [->] (u) edge (sibling);
            \path [->] (u) edge[bend left=20] (a);
            \path [->] (u) edge[bend left=30] (b);
            \path [->] (u) edge[bend left=40] (c);
            \path [->] (u) edge[bend left=50] (d);
            \path [->] (u) edge[bend left=60] (e);
        \end{scope}
    
        \begin{scope}[shift={(7.5-2,0)}]
            \draw[very thick, draw=lightgray] (-0.4,-0.5) -- (0.4,-0.5) -- (0,1) -- (-0.4,-0.5);

            \draw[very thick, draw=lightgray] (-0.4+1,-0.5) -- (0.4+1,-0.5) -- (0+1,1) -- (-0.4+1,-0.5);
            
            \node (b) at (1-0.15,0.6) {};
            \node (c) at (1-0.15,0.3) {};
            \node (d) at (1-0.15,0) {};
            \node (e) at (1-0.15,-0.3) {};

            \node[draw, shape = circle, fill = black, minimum size = 0.1cm, inner sep=1.6pt] (v1) at (0,1) {};
            
            \node[draw, shape = circle, fill = black, minimum size = 0.1cm, inner sep=1.6pt] (v2) at (1,1) {};
            
            \node (text cousin) at (2.6,2.4) {$l(u^{*})$};

            \draw[To-, black, dashed, thick] (v2)+(1.5mm,1mm) to [out=30,in=230] (text cousin);
        
            \node[draw, shape = circle, fill = black, minimum size = 0.1cm, inner sep=1.6pt] (w1) at (0.5,2) {};

            \node[draw, shape = circle, fill = black, minimum size = 0.1cm, inner sep=1.6pt] (w2) at (1.5,2) {};

            \node[draw, shape = circle, fill = black, minimum size = 0.1cm, inner sep=1.6pt] (w3) at (2,1) {};
            \node (text parent) at (2.55,1) {$u^{*}$};
            
            \node[draw, shape = circle, fill = black, minimum size = 0.1cm, inner sep=1.6pt] (root) at (1,3) {};
        
            \node(index) at (1,-1.5) {\textbf{G3}};
        \end{scope}
        
        \begin{scope}[shift={(7.5-2,0)}]
            \node[draw, shape = circle, fill = black, minimum size = 0.1cm, inner sep=1.6pt] (u) at (2.5,0) {};
            \node (text) at (3,0) {$u$};
        \end{scope}
        
        \begin{scope}[every edge/.style={draw, thick}]
            \path [-] (root) edge (u);
            \path [-] (root) edge (v1);
            \path [-] (w1) edge (v2);
        \end{scope}
        
        \begin{scope}[every edge/.style={draw=theme,very thick}]
            \path [->] (u) edge (v2);
            \path [->] (u) edge (b);
            \path [->] (u) edge (c);
            \path [->] (u) edge (d);
            \path [->] (u) edge (e);
        \end{scope}

        \begin{scope}[shift={(12,0)}]


        
            \draw[very thick, draw=lightgray] (-0.7-1.7,-0.5) -- (0.7-1.7,-0.5) -- (-1.7,2) -- (-0.7-1.7,-0.5);

            \draw[very thick, draw=lightgray] (-0.7,-0.5) -- (0.7,-0.5) -- (0,2) -- (-0.7,-0.5);

            \draw[very thick, draw=lightgray] (-0.7+1.7,-0.5) -- (0.7+1.7,-0.5) -- (1.7,2) -- (-0.7+1.7,-0.5);
            
            \node (cleft) at (-0.1-1.7,0.6) {};
            \node (dleft) at (-0.1-1.7,0.3) {};
            \node (eleft) at (-0.1-1.7,0) {};

            \node (cright) at (0.1+1.7,0.6) {};
            \node (dright) at (0.1+1.7,0.3) {};
            \node (eright) at (0.1+1.7,0) {};

            \node (bmid) at (-0.3,0.6) {};
            \node (cmid) at (0.3,0.6) {};
            \node (dmid) at (-0.5,0.3) {};
            \node (emid) at (0.5,0.3) {};
        
            \node(index) at (0,-1.5) {\textbf{G4}};
        \end{scope}
        
        \begin{scope}[shift={(12,0)}][draw, shape = circle, fill = black, minimum size = 0.1cm, inner sep=1.6pt]
            \node[draw, shape = circle, fill = black, minimum size = 0.1cm, inner sep=1.6pt] (u) at (0,-0.5) {};
            \node[draw, shape = circle, fill = black, minimum size = 0.1cm, inner sep=1.6pt] (u') at (0,2) {};
            \node[draw, shape = circle, fill = black, minimum size = 0.1cm, inner sep=1.6pt] (left cousin) at (-1.7,2) {};
            \node[draw, shape = circle, fill = black, minimum size = 0.1cm, inner sep=1.6pt] (right cousin) at (1.7,2) {};
            
            \node (text) at (0,-0.9) {$u$};
            \node (text middle) at (0,2.5) {$u'$};
            \node (text left) at (-1.7,2.5) {$l(u')$};
            \node (text right) at (1.7,2.5) {$r(u')$};
        \end{scope}
        
        \begin{scope}[every edge/.style={draw=theme,very thick}]
            \path [->] (u) edge[bend left=30] (left cousin);
            \path [->] (u) edge[bend left=15] (cleft);
            \path [->] (u) edge[bend left=15] (dleft);
            \path [->] (u) edge[bend left=15] (eleft);

            \path [->] (u) edge[bend right=30] (right cousin);
            \path [->] (u) edge[bend right=15] (cright);
            \path [->] (u) edge[bend right=15] (dright);
            \path [->] (u) edge[bend right=15] (eright);

            \path [->] (u) edge (u');
            \path [->] (u) edge (bmid);
            \path [->] (u) edge (cmid);
            \path [->] (u) edge (dmid);
            \path [->] (u) edge (emid);
        \end{scope}
        
    \end{tikzpicture}
    
    \caption{Types of edges in $\overrightarrow{G}^{h}_T$.}
    \label{fig: edges in GTr}
\end{figure}

\begin{definition}\label{newdef}
For a DFS-ordered tree $T$ and an integer $h\geq0$, define the directed graph \define{$\overrightarrow{G}^{h}_T$} on the vertex set $V(T)$ with the following arcs (see \Cref{fig: edges in GTr}). For every vertex $u$ of $T$, add to $\overrightarrow{G}^{h}_T$ the arcs:
\begin{enumerate}[label=\upshape \textbf{G\arabic{enumi}}, ref=\upshape \textbf{G\arabic{enumi}}]
     \item \label{G1}  $\overrightarrow{uw}$ for all $w\in D_{T}(u)$.
     \item \label{G2} $\overrightarrow{uw}$ for all left-siblings $w$ of $u$, and $\overrightarrow{uw'}$ for all $w'\in D_{T}(w)$.
     \item \label{G3} $\overrightarrow{uw}$ for the nearest-left cousin $w=l(u^*)$ of the parent $u^*$ of $u$, and $\overrightarrow{uw'}$ for all $w'\in D_{T}(w)$.
     \item \label{G4} $\overrightarrow{uu'}$ for the $h$-th ancestor $u'$ of $u$, $\overrightarrow{uw}$ for each $w\in\{l(u'),r(u')\}$, and $\overrightarrow{uw}$ for all $w\in D_{T}(l(u'))\cup D_{T}(u')\cup D_{T}(r(u'))$ such that $L(w)\leq L(u)$.
\end{enumerate}
\end{definition}

The graph \define{$h$-generated} by $T$, denoted by \define{$G^{h}_T$}, is the underlying undirected simple graph of $\overrightarrow{G}^{h}_T$. For $h=0$, we omit the superscript $h$ and say that $G_T$ is generated by $T$ (in which case, rule \ref{G4} adds no edges).

An induced subgraph $G^h_T[S]\subseteq G^h_T$ is \define{$T$-admissible} if $S$ is a $T$-admissible set. An embedding $\phi: V(T')\rightarrow V(G^{h}_T$) of a tree $T'$ is \define{$T$-admissible} if $V(G^{h}_T)\setminus\phi(V(T'))$ is a $T$-admissible set.

A \define{component} of a graph $G$ is the vertex set of a maximal connected subgraph of $G$. For a forest $T$ and a vertex $w\in V(T)$, a \define{$w$-component} of $T$ is a component of $T-w$. We make use of the following lemma about $w$-components. For a collection of sets $\{C_1,\dots, C_t\}$ and any subset $I$ of $[t]$, define $\define{$C_{I}$}:=\bigcup_{i\in I} C_i$.

\begin{lemma} \label{lem2}
For all integers $x\geq1$, for every forest $T$ with $|T|\geq x+1$, for every vertex $u\in V(T)$, there exists a vertex $w\in V(T)$ and a set $\{C_1,\dots, C_t\}$ of $w$-components of $T$ such that $u\notin C_{[t]}$ and $x\leq|C_{[t]}|\leq 2x-1$.
\end{lemma}
\begin{proof}
By adding more edges if necessary, we can assume that $T$ is a tree. Let $w_{0}:=u$. Let $C_{1},\ldots,C_{p}$ be the $w_{0}$-components of $T$. If $|C_{i}|\leq2x-1$ for all $i$, then it is easy to see that one can pick a desired collection of $w_{0}$-components from $\{C_{1},\ldots,C_{p}\}$. Otherwise, assume that $|C_{1}|\geq2x$. Let $w_{1}$ be the unique neighbour of $w_{0}$ contained in $C_{1}$. Let $D_{1},\ldots,D_{q}$ be the $w_{1}$-components of $T$ that do not contain $w_{0}$. If $|D_{i}|\leq2x-1$ for all $i$, then one can pick a desired collection of $w_{1}$-components from $\{D_{1},\ldots,D_{q}\}$. Otherwise, assume that $|D_{1}|\geq2x$. Let $w_{2}$ be the unique neighbour of $w_{1}$ contained in $D_{1}$. Since $T$ is a finite graph, repeating this process gives a sequence of vertices $w_{2},w_{3},\ldots,w_{\ell}$, such that each of the $w_{\ell}$-components of $T$ not containing $w_{\ell-1}$ have at most $2x-1$ vertices. Let $E_{1},\ldots,E_{r}$ be the $w_{\ell}$-components of $T$ that do not contain $w_{\ell-1}$. Then one can pick a desired collection of $w_{\ell}$-components from $\{E_{1},\ldots,E_{r}\}$.
\end{proof}

A collection $\{C_{1},\dots, C_{t}\}$ of $w$-components of a forest $T$ is \define{$(x,y)$-feasible} if $x\leq |C_{[t]}\cup \{w\} |\leq x+y-2$. The collection is \define{$(x,y)$-critical} if $t\geq 2$, $x+y-2 \leq  |C_{[t]}| \leq  2x-3$, and $|C_{I}|\leq  x-2$ for all $I\subsetneq [t]$. Note that an $(x,y)$-critical collection satisfies $|C_i|\geq y$ for all $i\in [t]$.  Furthermore, \Cref{lem2} ensures that if $|T|\geq x+1$, then an $(x+1,x+1)$-feasible collection always exists in $T$ whenever $x\geq 1$. The following lemma guarantees the existence of either an $(x,y)$-feasible collection or an $(x,y)$-critical collection.
\begin{lemma} \label{lem: feasible critical}
   For all integers $x>y\geq2$, for any forest $T$ with $|T|\geq x$ and a vertex $u\in V(T)$, there exists a vertex $w\in V(T)$ and a set $\{C_1,\dots, C_t\}$ of $w$-components of $T$ such that $u\notin C_{[t]}$ and $\{C_1,\dots, C_t\}$ is either $(x,y)$-feasible or $(x,y)$-critical.
\end{lemma}
\begin{proof}
Since $|T|\geq x$, by \Cref{lem2}, there exists a vertex $w_{0}\in V(T)$ and a set $\{C_1,\dots, C_t\}$ of $w_{0}$-components of $T$ such that $u\notin C_{[t]}$ and $x-1\leq|C_{[t]}|\leq 2x-3$. 

We are done if $\{C_1,\dots, C_t\}$ is $(x,y)$-feasible, so assume that $x+y-2\leq|C_{[t]}|\leq 2x-3$. 
 Let $\{D_{1},\ldots,D_{q}\}$ be a minimal subset of $C_{[t]}$ such that $|D_{[q]}|\geq x+y-2$. Note that $|D_{[q]}|\leq|C_{[t]}|\leq 2x-3$. Assume $q\geq2$. 
 If there exists $I\subsetneq[q]$ such that $|D_{I}|\geq x-1$, then by the minimality of $D_{[q]}$, $x\leq |D_{I}\cup\{w_{0}\}|\leq x+y-2$, implying $\{D_{i}\}_{i\in I}$ is $(x,y)$-feasible.
 Otherwise, $D_{[q]}$ is $(x,y)$-critical by definition. 
 
Now assume $q=1$. Let $w_{1}$ be the unique neighbour of $w_{0}$ in $D_{1}$. Let $E_{1},\ldots,E_{r}$ be the $w_{1}$-components of $T[D_{1}]$. Since $E_{[r]}=D_{1}\setminus\{w_{1}\}$ and $y\geq2$, $u\not\in E_{[r]}$ and $x-1\leq|E_{[r]}|\leq2x-3$. Since $T$ is a finite graph, repeating this process with $\{E_1,\dots, E_r\}$ playing the role of $\{C_1,\dots, C_t\}$
gives a sequence of vertices $w_{1},w_{2},\ldots,w_{\ell}$, such that some set of $w_{\ell}$-components of $T$ not containing $u$ is either $(x,y)$-feasible or $(x,y)$-critical as desired.
\end{proof}

\begin{definition}\label{def: treewidth}
   A \define{tree-decomposition} of a graph $G$ is a collection of sets $(B_x \subseteq V(G) : x\in V(T))$ (called \define{bags}) indexed by the vertices of a tree $T$, such that each of the following hold: 
\begin{enumerate}[label=\upshape \textbf{D\arabic{enumi}}, ref=\upshape \textbf{D\arabic{enumi}}]
   \item \label{D1} for all $u\in V(G)$, $T_u:=T[\{x\in V(T):u\in B_x\}]$ is a nonempty subtree of $T$, and
   \item \label{D2} for all $uv\in E(G)$, there exists $x\in V(T)$ such that $\{u,v\}\subseteq B_x$.
\end{enumerate}
The \define{width} of the tree-decomposition $(B_x : x\in V(T))$ equals $\max_{x\in V(T)} |B_x|-1$. The \define{treewidth} of $G$, denoted by \define{$\tw(G)$}, is the minimum width of a tree-decomposition of $G$.
\end{definition}

We use the following theorem of \citet{wood} to prove \Cref{thm: treewidth}.
\begin{theorem}\label{eq: tree decomp}
Every graph $H$ with ${\rm tw}(H)\leq k$ has a tree-decomposition $(B_x:x\in V(T))$ with $|T|\leq \frac{1}{k}|H| +1 $ and $ |B_x|\leq 3k$ for all $x\in V(T)$.
\end{theorem}



\section{The Chung--Graham Error}\label{sec: 3}

This section explains the error by \citet{grahamthird}, and provides a fix that implies the bound in \eqref{eqn: corrected bound}. Our fix is similar but not the same as that of \citet{frati}. We first introduce some definitions from \cite{grahamthird}. Let \define{$B(k)$} denote the (rooted) perfect binary tree with level $k$ that is DFS-ordered.


\begin{definition}\label{olddef}
For any integer $k\geq0$, define the graph \define{$G(k)$} on the vertex set of $B(k)$ by adding the following edges. For each $u\in V(B(k))$, add the edges:
\begin{enumerate}[label=\upshape\textbf{(\roman*)}, ref=\upshape\textbf{(\roman*)}]
     \item\label{olddef: i} $uw$ for all $w\in D_{B(k)}(u)$,
     \item\label{olddef: ii} $uw$ for all left-siblings $w$ of $u$, and $uw'$ for all $w'\in D_{B(k)}(w)$, and
     \item\label{olddef: iii} $uw$ for the left-sibling $w$ of the parent $u^{*}$ of $u$ (if they exist), and $uw'$ for all $w'\in D_{B(k)}(w)$.
\end{enumerate} 
\end{definition}




Note that the graph $G(k)$ of \Cref{olddef} is not the same as the graph $G_{B(k)}$ of \Cref{newdef}. In particular, $w$ in \labelcref{olddef: iii} of \Cref{olddef} is the left-sibling of $u^{*}$, instead of the nearest-left cousin of $u^{*}$ like in \labelcref{G3} of \Cref{newdef}.

A subgraph $G$ of $G(k)$ that is induced by a $B(k)$-admissible set is called \define{$G(k)$-admissible}. Furthermore, for any vertex $v\in V(B(k))$, define
$\define{$D_{G}(v)$}:=V(G)\cap D_{B(k)}(v)$. The following is the main lemma in \cite{grahamthird}.

\begin{lemma}[labelled $(*)$ in {\cite[page 205]{grahamthird}}] 
\label{oldlemma}
For any $B(k)$-admissible graph $G$, for any tree $T$ with $|T|\leq|G|$, for any vertex $v\in V(T)$, there exists an embedding $\lambda : V(T)\to V(G)$ such that $V(G)\setminus\lambda(V(T))$ is $B(k)$-admissible, and the location of $\lambda(v)$ can be determined as follows:
\begin{enumerate}
    \item[\scriptsize$\bullet$] if $|T|=|G|$, then $\lambda(v)$ is the root of $B(k)$, and
    
    \item[\scriptsize$\bullet$] if $|T|<|G|$, then $\lambda(v)$ is the last vertex in $\lambda(V(T))$ (as per the DFS order) such that
    \begin{align*}
        |D_{G}(\lambda(v))| < |T|\leq  |D_{G}(\lambda(v)^*)|\;,
    \end{align*}
where the last inequality need only hold when $\lambda(v)^*$ exists.
\end{enumerate}
\end{lemma}

    

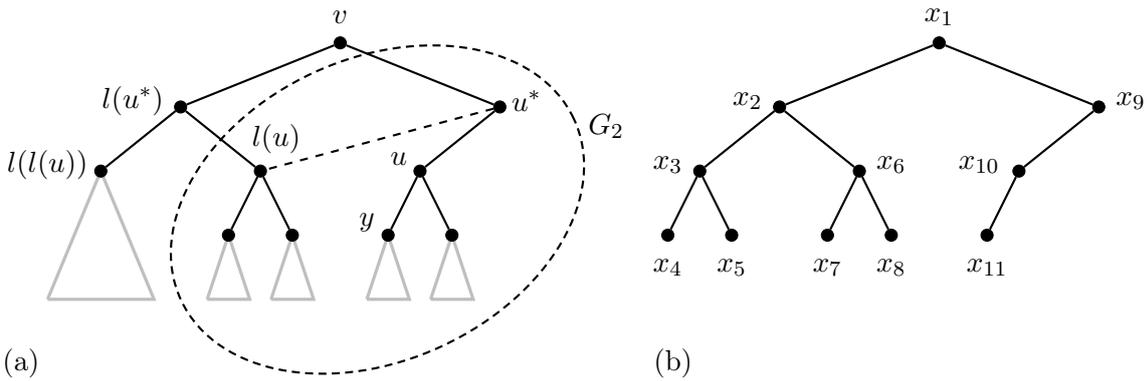
\begin{figure}[!ht]
%
\centering
    \begin{tikzpicture}[/tikz/xscale=0.7, /tikz/yscale=0.85]
    \begin{scope}
        \draw[very thick, draw=lightgray] (-1-4.5,-1) -- (1-4.5,-1) -- (-4.5,1) -- (-1-4.5,-1);
        \draw[very thick, draw=lightgray] (-0.4-1.5-0.6,-1) -- (0.4-1.5-0.6,-1) -- (-1.5-0.6,0) -- (-0.4-1.5-0.6,-1);
        \draw[very thick, draw=lightgray] (-0.4-1.5+0.6,-1) -- (0.4-1.5+0.6,-1) -- (-1.5+0.6,0) -- (-0.4-1.5+0.6,-1);
        \draw[very thick, draw=lightgray] (-0.4+1.5-0.6,-1) -- (0.4+1.5-0.6,-1) -- (1.5-0.6,0) -- (-0.4+1.5-0.6,-1);
        \draw[very thick, draw=lightgray] (-0.4+1.5+0.6,-1) -- (0.4+1.5+0.6,-1) -- (1.5+0.6,0) -- (-0.4+1.5+0.6,-1);
    \end{scope}
    
    \begin{scope}[every node/.style={draw, shape = circle, fill = black, minimum size = 0.1cm, inner sep=1.6pt}]
        \node (v) at (0,3) {};
        \node (lustar) at (-3,2) {};
        \node (ustar) at (3,2) {};
        \node (llu) at (-4.5,1) {};
        \node (lu) at (-1.5,1) {};
        \node (u) at (1.5,1) {};
        
        \node (w1) at (-1.5-0.6,0) {};
        \node (w2) at (-1.5+0.6,0) {};
        \node (w3) at (1.5-0.6,0) {};
        \node (w4) at (1.5+0.6,0) {};
    \end{scope}
    
    \begin{scope}[every edge/.style={draw, thick}]
        \path [-] (v) edge (lustar);
        \path [-] (v) edge (ustar);
        \path [-] (lustar) edge (llu);
        \path [-] (lustar) edge (lu);
        \path [-] (lu) edge (w1);
        \path [-] (lu) edge (w2);
        \path [-] (ustar) edge (u);
        \path [-] (u) edge (w3);
        \path [-] (u) edge (w4);

        \path [thick, draw=black, dashed] (ustar) edge (lu);
    \end{scope}

    \begin{scope}[every node/.style={ellipse, thick, draw, black}]
        \node[densely dashed, minimum width=5.7cm, minimum height=4cm, rotate=25] (G2) at (0.7,0.4) {};
    \end{scope}

    \begin{scope}
        \node (v) at (0,3+0.4) {$v$};
        \node (lustar) at (-3-0.9,2+0.1) {$l(u^{*})$};
        \node (ustar) at (3+0.5,2+0.1) {$u^{*}$};
        \node (llu) at (-4.5-1,1+0.1) {$l(l(u))$};
        \node (lu) at (-1.5+0.3,1+0.5) {$l(u)$};
        \node (u) at (1.5-0.4,1+0.2) {$u$};
        
        \node (w1) at (-1.5-0.6,0) {};
        \node (w2) at (-1.5+0.6,0) {};
        \node (w3) at (1.5-0.6-0.4,0.2) {$y$};
        \node (w4) at (1.5+0.6,0) {};

        \node (G2) at (5,2-0.3) {$G_{2}$};
        \node (a) at (-6,-2) {(a)};

    \end{scope}
        
    \end{tikzpicture}
    \begin{tikzpicture}[/tikz/xscale=0.7, /tikz/yscale=0.85]
    \begin{scope}
        \node[densely dashed, minimum width=5.7cm, minimum height=4cm, rotate=25, ellipse, thick, draw, white] (G2) at (0.7,0.4) {};
    \end{scope}
    
    \begin{scope}[every node/.style={draw, shape = circle, fill = black, minimum size = 0.1cm, inner sep=1.6pt}]
        \node (x1) at (0,3) {};
        \node (x2) at (-3,2) {};
        \node (x3) at (-4.5,1) {};
        \node (x4) at (-4.5-0.6,0) {};
        \node (x5) at (-4.5+0.6,0) {};
        \node (x6) at (-1.5,1) {};
        \node (x7) at (-1.5-0.6,0) {};
        \node (x8) at (-1.5+0.6,0) {};
        \node (x9) at (3,2) {};
        \node (x10) at (1.5,1) {};
        \node (x11) at (1.5-0.6,0) {};
    \end{scope}
    
    \begin{scope}[every edge/.style={draw, thick}]
        \path [-] (x1) edge (x2);
        \path [-] (x1) edge (x9);
        \path [-] (x2) edge (x3);
        \path [-] (x2) edge (x6);
        \path [-] (x3) edge (x4);
        \path [-] (x3) edge (x5);
        \path [-] (x6) edge (x7);
        \path [-] (x6) edge (x8);
        \path [-] (x9) edge (x10);
        \path [-] (x10) edge (x11);
    \end{scope}

    \begin{scope}
        \node (x1) at (0,3+0.4) {$x_{1}$};
        \node (x2) at (-3-0.6,2+0.1) {$x_{2}$};
        \node (x3) at (-4.5-0.6,1+0.1) {$x_{3}$};
        \node (x4) at (-4.5-0.6,0-0.5) {$x_{4}$};
        \node (x5) at (-4.5+0.6,0-0.5) {$x_{5}$};
        \node (x6) at (-1.5+0.6,1+0.1) {$x_{6}$};
        \node (x7) at (-1.5-0.6,0-0.5) {$x_{7}$};
        \node (x8) at (-1.5+0.6,0-0.5) {$x_{8}$};
        \node (x9) at (3+0.6,2+0.1) {$x_{9}$};
        \node (x10) at (1.5-0.75,1+0.1) {$x_{10}$};
        \node (x11) at (1.5-0.6,0-0.5) {$x_{11}$};
        \node (b) at (-5,-2) {(b)};
    \end{scope}
        
    \end{tikzpicture}
        \caption{(a) Defining $G_{2}$. 
(b) The $11$-vertex $G(3)$-admissible graph.}
    \label{fig: counterexample}
\end{figure}

The error in the proof of \Cref{oldlemma} occurs in Part II, Case (ii), Subcase (b), as illustrated in \Cref{fig: counterexample}(a) (which is the same as Figure~7 in \cite[page 207]{grahamthird}). \citet{grahamthird} claimed that the graph $G_2$ in \Cref{fig: counterexample}(a) is isomorphic to a $G(k)$-admissible graph. However, this claim is false. We now construct a small counterexample. Consider the $G(3)$-admissible graph on $11$ vertices with the DFS preorder traversal $x_1,\dots, x_{11}$ as shown in \Cref{fig: counterexample}(b). In this case, $G_{2}$ corresponds to the subgraph induced by the last six vertices $x_6,\dots, x_{11}$. Since $x_{6}$ is not the left-sibling of $x_{10}$ (the parent of $x_{11}$), \ref{olddef: iii} does not add the edges $x_6x_{11}, x_7x_{11}$ and $x_8x_{11}$. Thus, $G_{2}$ is isomorphic to a non-complete graph.
It is routine to check that for any $\ell$, every $6$-vertex $G(\ell)$-admissible graph is isomorphic to the complete graph $K_{6}$, implying $G_{2}$ is not admissible in the given example.


Nonetheless, the proof can be corrected. By replacing $G(k)$ with the graph $G_{B(k)}$ (defined in \Cref{newdef}), the induction argument for \Cref{oldlemma} works without additional modifications. Therefore, it can be shown that every $n$-vertex $G_{B(k)}$-admissible graph is $\mathcal{T}_{n}$-universal as desired. Moreover, one may check that this correction increases the number of edges in the resulting $n$-vertex universal graph by at most $n k + \bigo(n)$.  Hence, \labelcref{eqn: corrected bound} holds.

\section{Universality for Treewidth-$k$ Graphs}\label{sec: universality for treewidth}

In this section we prove \Cref{thm: treewidth} for treewidth-$k$ graphs. We first show the lower bound.

\begin{proof}[Proof of lower bound in \Cref{thm: treewidth}] Let $c:=e^{17/4}-1$. Let $U$ be a $\mathcal{T}_{n}^{k}$-universal graph. We first assume that $1<n/k\leq c$. Note that $U$ contains a copy of the complete graph $K_{k+1}$. So $e(U)\geq{k+1\choose2}=\Omega(k^{2})$. Since $k\geq n/c$ and $1\geq\frac{\ln(n/k+1)}{\ln(c+1)}$, $e(U)=\Omega(k\cdot\frac{n}{c}\cdot\frac{\ln(n/k+1)}{\ln(c+1)})=\Omega(kn\log(n/k+1))$.

Now assume that $n/k\geq c$. Note that $U$ contains the disjoint union of $j$ copies of the complete bipartite graph $K_{k,\lfloor\frac{n}{j}\rfloor-k}$ for each $j\leq \lfloor n/k\rfloor$. Thus, there are at least $jk$ vertices of degree at least $\lfloor\frac{n}{j}\rfloor-k$ in $U$ for each $j\leq \lfloor n/k\rfloor$. So the number of edges in $U$ is at least
\[\frac{1}{2}\sum_{j=1}^{\lfloor\frac{n}{k}\rfloor}k\Big(\Big\lfloor\frac{n}{j}\Big\rfloor-k\Big)\geq\frac{1}{2}kn\sum_{j=1}^{\lfloor\frac{n}{k}\rfloor}\frac{1}{j}-kn\geq\frac{1}{2}kn\ln(\lfloor n/k\rfloor+2)-\frac{1}{2}\cdot\frac{kn}{\lfloor\frac{n}{k}\rfloor+1}-\frac{1}{2}\cdot\frac{kn}{\lfloor\frac{n}{k}\rfloor+2}-kn\;.\]
Since $n/k>1$, $\frac{1}{2}\cdot\frac{kn}{\lfloor\frac{n}{k}\rfloor+1}+\frac{1}{2}\cdot\frac{kn}{\lfloor\frac{n}{k}\rfloor+2}+kn\leq\frac{17}{12}kn$. Since $n/k\geq c$, $\frac{17}{12}kn\leq\frac{1}{3}kn\ln(n/k+1)$, implying $e(U)\geq\frac{1}{2}kn\ln(n/k+1)-\frac{1}{3}kn\ln(n/k+1)=\Omega(kn\log(n/k+1))$ as desired.
\end{proof}

To prove the upper bound in \Cref{thm: treewidth}, we begin by making some crucial observations about the proof of \Cref{oldlemma} in \cite{grahamthird} when applied to $G_{B(k)}$ as described in \Cref{sec: 3}. This proof inductively embeds a given tree $T$ on at most $n$ vertices into an $n$-vertex $B(k)$-admissible graph $G$. Several cases are considered in \cite{grahamthird}, but in general, the strategy is to find a vertex $u$ of $T$ and a collection $C_{1},\ldots,C_{t}$ of $u$-components of $T$ satisfying certain size conditions, and then embedding the forests $T_{1}:=T[\{u\}\cup C_{[t]}]$ and $T_{2}:=T-V(T_{1})$ into $G$ by induction. In some cases, $T_{1}$ might consist of a single vertex, or $T_{2}$ might be the empty tree. In each of the induction applications, the proof repeatedly divides the given tree into two.

Let $\lambda:V(T)\to V(G)$ be the final embedding obtained from \Cref{oldlemma}. Then we can define the recursion tree \define{$T_{\rm aux}$} (for $\lambda$) as follows. It is rooted at $u$, the left child of $u$ is the vertex $u_{1}$ found in $T_{1}$ (for the induction application in $T_{1}$), the right child of $u$ is the vertex $u_{2}$ found in $T_{2}$, and so on. Note that $V(T_{\rm aux})=V(T)$.

\begin{observation}\label{obs: Taux embedding}
    For any vertex $w$ of $T_{\rm aux}$, the proof of \Cref{oldlemma} guarantees that $\lambda(w)$ is adjacent to every vertex in $\lambda\left(D_{T_{\rm aux}}(w)\right)$ in $G$.
\end{observation}

The next observation follows from the definition of $T_{\rm aux}$.

\begin{observation}\label{obs: Taux ancestor}
    If two vertices $x,x'\in V(T_{\rm aux})$ have a closest common ancestor $y$ in $T_{\rm aux}$ such that $y\notin \{x,x'\}$, then $x$ and $x'$ belong to distinct $y$-components of $T$.
\end{observation}

Let \define{$T_{\rm aux}^{*}$} be the graph obtained from $T_{\rm aux}$ by adding edges from each vertex $x$ of $T_{\rm aux}$ to every vertex in $D_{T_{\rm aux}}(x)$. \Cref{obs: Taux embedding} implies that $\lambda$ is an embedding of $T_{\rm aux}^{*}$ into $G$. Furthermore, by \Cref{obs: Taux ancestor}, for any edge $xy\in E(T)$, one of $x,y$ is a $T_{\rm aux}$-ancestor of the other, implying $T$ is a subgraph of $T_{\rm aux}^*$. We apply these observations in the following argument.

\begin{proof}[Proof of upper bound in \Cref{thm: treewidth}]
Let $m:=\lfloor\frac{n}{k}+1\rfloor$. Since $n>k$, $m\geq2$ and $k^{2}m\leq2kn$. Consider the $m$-vertex $B(\lceil\log_{2}m\rceil)$-admissible graph $G$. Then $G$ is $\mathcal{T}_{m}$-universal by \Cref{oldlemma}, and $e(G)\leq cm\log_{2} m$ for some constant $c>0$. Let $G'$ be the graph obtained from $G$ by replacing each vertex $x$ of $G$ with a clique $C_{x}$ of size $3k$. Observe that each edge of $G$ becomes a copy of the complete bipartite graph $K_{3k,3k}$ in $G'$. Therefore, $e(G')=9k^{2}e(G)+{3k\choose2}|G|\leq 9ck^2m\log_{2}m+\frac{9}{2}k^{2}m\leq(9c+\frac{9}{2})k^{2}m\log_{2}m$, where the last inequality holds since $m\geq2$. Since $k^{2}m\leq2kn$, $e(G')=\bigo(kn\log(n/k+1))$ as desired.

We claim that $G'$ is $\mathcal{T}_{n}^{k}$-universal. Let $H$ be an $n$-vertex treewidth-$k$ graph. By \Cref{eq: tree decomp}, $H$ has a tree-decomposition $(B_{x}:x\in V(T))$ with width at most $3k-1$ and $|T|\leq m$. Let $\lambda:V(T)\to V(G)$ be the embedding from \Cref{oldlemma}. Define $T_{\rm aux}$ and $T_{\rm aux}^{*}$ (for this $\lambda$) as above. For each vertex $v\in V(H)$, define $f(v)$ to be the unique vertex of $T_{v}=T[\{x\in V(T):v\in B_{x}\}]$ with minimum level in $T_{\rm aux}$; $f(v)$ is well-defined, since if $v\in B_{x}\cap B_{x'}$ and $L_{T_{\rm aux}}(x)=L_{T_{\rm aux}}(x')$, then \Cref{obs: Taux ancestor} and \labelcref{D1} in \Cref{def: treewidth} together imply that the closest common ancestor of $x$ and $x'$ in $T_{\rm aux}$ is contained in $T_{v}$.

Since each bag $B_{x}$ of $T$ contains at most $3k$ vertices, there is an injection $\pi:V(H)\to V(G')$ that sends each vertex $v$ of $V(H)$ to a vertex in $C_{\lambda(f(v))}$. It remains to show that $\pi$ is a homomorphism. Let $uv\in E(H)$. Let $y$ be the closest common ancestor of $f(u)$ and $f(v)$ in $T_{\rm aux}$. Then $B_{y}$ separates $u$ and $v$ in $H$. Since $uv\in E(H)$, it must be that $u$ or $v$ is in $B_{y}$, implying $f(u)$ or $f(v)$ equals $y$. Therefore, $f(u)f(v)\in E(T_{\rm aux}^{*})$. Since $\lambda$ is an embedding of $T_{\rm aux}^{*}$ into $G$, it follows that $\pi(u)\pi(v)\in E(G')$. Hence, $\pi$ embeds $H$ into $G'$ as desired.
\end{proof}

Note that the graph $G'$ above has $3km\leq3(n+k)\leq6n$ vertices. Moreover, any universal graph for $n$-vertex treewidth-$k$ graphs that has $\bigo(n)$ vertices, has treewidth $\Omega(k\log(n/k+1))$ whenever $n>k$. This lower bound follows immediately from the fact that every $n$-vertex graph of treewidth $t$ has at most $tn - {t+1 \choose 2}$ edges. On the other hand, we now show that the treewidth of the universal graph $G'$ is at most $51k\log_{2}(n/k+1)$.

\begin{proof} For this upper bound, we first construct a tree-decomposition of $G_{B(\ell)}$ indexed by the path $T=(1,2,\dots,2^{\ell})$, where $\ell := \lceil \log_2 m \rceil$. Note that $B(\ell)$ 
has exactly $2^{\ell}$ leaves $x_1, \dots, x_{2^{\ell}}$ ordered from left to right.
For each $i \in [2^{\ell}]$, let $B_i = N_{G_{B(\ell)}}(x_i) \cup \{x_i\}$ be the set consisting of $x_i$ and all its neighbours. We claim $(B_1,B_2,\dots,B_{2^\ell})$ is a tree-decomposition of $G_{B(\ell)}$. For each level $d\in\{0,\ldots,\ell\}$, consider the unique ancestor $y$ of $x_i$ such that $y$ has level $d$. Then the vertices at level $d$ in $B_i$ are the children of $y^*$ or the children of the nearest-right cousin of $y^*$ if $d<\ell-1$. If $d\in \{\ell-1,\ell\}$, then vertices at level $d$ in $B_i$ may additionally be the children of the nearest left-cousin of $y^*$ or left siblings of $y$. 
Hence, at most four vertices at level $d$ belong to $B_i$ if $d<\ell-1$ and at most six vertices at level $d$ belong to $B_i$ if $d\in \{\ell-1,\ell\}$ so such a tree-decomposition has width at most $4\ell + 8$. Indeed, for any vertex $u$ in $G_{B(\ell)}$, its neighbours in $\{x_1, \dots, x_{2^{\ell}}\}$ are either descendants of $u$ by \ref{G1}, descendants of left siblings of $u$ by \ref{G2}, descendants of the nearest left-cousin of $u^*$, children of $r(u)$ (if $u$ has level $\ell-1$), and right siblings of $u$ (if $u$ has level $\ell$) by \ref{G3}. These correspond to consecutive vertices $\{x_i, x_{i+1}, \dots, x_{i+j}\}$ for some $i$ and $j$, which form a subpath of $T$. Thus \ref{D1} holds.
Moreover, if $\overrightarrow{uw}$ is an arc in $\overrightarrow{G}_{B(\ell)}$, then $\overrightarrow{uw'}$ is also an arc for every descendant $w'$ of $w$. Hence, there exists a leaf-descendant $x_i$ of $w$ adjacent to both $u$ and $w$ in $G_{B(\ell)}$, verifying that \ref{D2} holds, and thus $(B_i : i \in V(T))$ is indeed a valid tree-decomposition.
Once we have this tree-decomposition, blowing up each vertex into a clique of size $3k$ naturally yields another tree-decomposition of the blown-up graph $G''$. The new decomposition has width at most $3k(4\ell+9)-1 < 12k \log_{2} m+39k\leq51k\log_{2}m$, where the last inequality holds because $m\geq2$. Since the graph $G'$ in the above proof is a subgraph of $G''$, it follows that $G'$ has treewidth at most $51k\log_{2}m\leq51k\log_{2}(n/k+1)$, as desired.
\end{proof}

This completes the proof of \Cref{thm: treewidth}.

\section{Universality for Trees}\label{sec:4}

\subsection{Idea and Outline of the Proof}\label{sec: overview}
 
 As explained in \Cref{sec: 3}, \citet{grahamthird} considered a graph generated by a perfect binary tree $B(k)$, leading (after applying our fix) to a universal graph with $\frac{7}{2}n\log_2 n + \bigo(n)$ edges. A natural direction for further improvement is to replace $B(k)$ with more general underlying tree $T$. There are, however, two separate aspects to this approach. First, we need to understand which choices of $T$ still guarantee universality. Second, among those choices, we would like to identify trees that lead to fewer edges in the resulting graph. In what follows, we examine these two aspects in turn. 

 To reiterate the rough idea in \cite{grahamthird}, suppose that we want to embed a tree $T'$ into $G_T$, where $G_T$ looks like \Cref{fig: counterexample}(a). Then we can use \Cref{lem2} to find a vertex $w$ and some $w$-components $C_1,\dots, C_t$ of $T'$ such that $C_{[t]}$ has a size between $\nu_T(u)$ and $2\nu_T(u)$. 
Take a subgraph $G_2$ induced by $D_{T}[l(u)]\cup D_{T}[u]\cup \{u^*\}$ as in \Cref{fig: counterexample}(a). By showing this is isomorphic to a $T$-admissible graph, by induction we can find a $T$-admissible embedding $\phi$ of $T'_0 = T'[\{w\}\cup C_{[t]}]$ into $G_2$ as desired. Furthermore, we ensure that $w$ maps to $u$. Once again, by induction we find an embedding $\phi'$ of $T'_1= T'- V(T'_0)$ into $G_{T} - \phi(V(T'_0))$. Since all the edges between $T'_0$ and $T'_1$ are incident to $w$, and $u=\phi(w)$ is adjacent to every vertex in $G_T$, we can concatenate $\phi$ and $\phi'$ to obtain a desired embedding. One crucial condition for this argument to work is that $T'_0$ has fewer than $|G_2|$ vertices. This condition is implied by $\nu_T(l(u))\geq \nu_T(u)$. Thus, this argument requires $l(u)$ to have at least as many descendants as $u$ for every vertex $u$.

 On the other hand, $\nu_T(l(u))$ being smaller than $\nu_T(u)$ helps to obtain a construction with fewer edges. Roughly speaking, this is because the number of edges we obtain from \Cref{newdef}~\ref{G2} and \ref{G3} becomes smaller if the left-siblings (or the nearest left-cousin) of $u$ has a smaller number of descendants. 
If we allow the ratio between $\nu_T(l(u))$ and $\nu_T(u)$ to be as large as $K=\Omega(\varepsilon^{-3})$, then we can construct a tree $T$ with $e(G_T)\leq (\frac{5}{2}+\varepsilon) n \log_2 n$. However, we were not able to prove that such a graph contains all possible spanning trees, as the imbalance imposes additional difficulty in embedding an arbitrary tree $T'$ into $G^h_T$.

 The crucial point is that, in the induction process, $C_{[t]}\cup \{w\}$ might be too large to embed into the union of two trees induced by $D_T[l(u)]$ and $D_T[u]$. 
So, we sometimes have to merge three trees $D_T[l(l(u))], D_T[l(u)], D_T[u]$ (instead of two) into one tree to invoke the induction hypothesis.
However, as some edges between $D_T[u]$ and $D_T[l(l(u))]$ might not be present, embedding the tree becomes much more difficult, if not impossible.  
In fact, for large imbalance $K$, we are able to find a desired embedding only by adding additional edges to $G_T$. In addition to \Cref{newdef}, if we add edges from each vertex $u$ to its nearest-right cousin and all of its descendants in $G_T$, a technical analysis yields that such additions produce a universal graph. Moreover, we can construct such a tree $T$ with  $e(G_T)\leq (3+\varepsilon) n \log_2 n$. Although this is a significant improvement from the $\frac{7}{2}n\log_2 n+\bigo(n)$ bound by Chung and Graham, we do not present this construction. Instead, we present a better construction.

In this section, we show how to obtain a better upper bound by taking $K=4$. The main issue for embedding the subtree $T'_0 = T'[\{w\}\cup C_{[t]}]$ into the union of three trees $D_T[l(l(u))], D_T[l(u)], D_T[u]$ is that the edges $ww_1, ww_2,\dots, w w_t$ between the vertex $w$ and the $w$-components $C_1,\dots, C_t$, respectively, are problematic if $w$ is mapped into $D_T[l(u)]$ or $D_{T}[l(l(u))]$ since some $w_i$ must be mapped into $D_T[u]$. In order to ensure that $\phi(w)$ and $\phi(w_i)$ are adjacent, we map each vertex $w_i$ such that $L_T(\phi(w))$ is small and $L_T(\phi(w_i))$ is less than $L_T(\phi(w))+2$. Then the additional condition \ref{G4} in \Cref{newdef} for $h=2$ will guarantee that $\phi(w)$ is adjacent to $\phi(w_i)$ as desired. Furthermore, we utilise \Cref{lem: feasible critical} to obtain a more sophisticated partition of $T'_0=T'[\{w\}\cup C_{[t]}]$ that is useful for obtaining a desired embedding.
In order to obtain such a partition, we crucially use the fact that the ratio between $\nu_T(l(u))$ and $\nu_T(u)$ is at most four.
These arguments are encapsulated in \Cref{strong}, which shows that one can embed a tree $T'$ into $G^2_T$ so that certain vertices $x_1,x_2$ are guaranteed to be embedded such that  $L_T(\phi(x_1))$ and $L_T(\phi(x_2))$ are small (provided that certain conditions are met). \Cref{strong} is then used in the construction of the desired universal graph with $\tfrac{14}{5}n\log_2 n + \bigo(n)$ edges. 

\subsection{\boldmath $(K,s)$-Trees}

In this section, we define a tree with a controlled imbalance between $\nu(u)$ and $\nu(l(u))$. Moreover, we prove that the graph $2$-generated by such an $n$-vertex tree is a universal graph for $\mathcal{T}_n$ in the next section.

\begin{definition}\label{def: Ks tree}
    A DFS-ordered tree $T$ is a \define{$(K,s)$-tree} if for every vertex $u\in V(T)$,
\begin{enumerate}[label=\upshape \textbf{T\arabic{enumi}}, ref=\upshape \textbf{T\arabic{enumi}}]
        \item \label{T1} 
        if $l(l(u))$ exists, then $\nu(l(l(u)))+\nu(l(u)) \geq \nu(u)$, 
        \item \label{T2}
        if $l(u)$ exists, then $\nu(l(u)) \geq \frac{1}{K} \nu(u)$,
        \item \label{T3} if $u$ is not the right-most (in the sense of $DFS$-order) vertex on its level, then for every vertex $u'\in V(T)$ with $L(u')\geq L(u)+s$, $\nu(u)\geq \nu(u')$, and
        \item \label{T4} if $u$ has a child, then $l(u)$ has a child.
    \end{enumerate}
\end{definition}
 
In other words, $T$ resembles a perfect binary tree in the sense that each vertex $u$ at level $\ell$ has roughly the same  number of descendants (up to some  multiplicative constant), and a vertex on a lower level has more descendants than a vertex on a higher level if the level difference is at least $s$. 
 Moreover, \ref{T1} ensures that if $\nu(l(u))$ is smaller than $\nu(u)$, then $\nu(l(l(u)))$ must be somewhat large, so that for any two consecutive vertices in level $\ell$ at least one of the $\nu$-values is large. We observe that the notion of $(K,s)$-tree is compatible with admissible sets.

\begin{observation}
    For any $(K,s)$-tree $T$, for any $T$-admissible set $U$, $T[U]$ is a $(K,s)$-tree.
\end{observation}

The following lemma gives another way to obtain a $(K,s)$-tree from a $(K,s)$-tree.

\begin{lemma} \label{lem: tree merge}
    For any $(K,s)$-tree $T$ and any set $\{u_{1},\ldots,u_{t}\}$ of vertices of $T$ such that for each $i$, $u_{i+1}=r(u_{i})$, let $Q=D_T[u_1]\cup \dots \cup D_T[u_{t}]$ and let $T^{*}$ be the new tree obtained from $T[Q]$ by adding a new root vertex $v$ adjacent to $u_1,\dots, u_t$. Then
    \begin{enumerate}[label=(\arabic*)]
        \item $T^{*}$ is a $(K,s)$-tree, and
        \item if $u_{1}^{*}=u_{t}^{*}$ or $u_{1}^{*}=l(u_{t}^{*})$, then for all $h\geq0$, $G^{h}_{T^*}$ is isomorphic to a subgraph of $G^{h}_T[Q \cup \{u_{t}^*\}]$, where $v$ is mapped to $u_{t}^*$.
    \end{enumerate}
\end{lemma}
\begin{proof}
    From left to right, the children of $v$ in $T^{*}$ are $u_{1},u_{2},\ldots,u_{t}$. Then DFS preorder traversal yields a linear order on $V(T^{*})$ in which $v$ comes first, followed by the remaining vertices in the same relative order as they appear in the linear order of $V(T)$.
    
    Observe that for every non-root vertex $u$ of $T^{*}$, $D_{T^{*}}[u]=D_{T}[u]$, $L_{T^{*}}(u)=L_{T}(u)-L_{T}(u^{*}_1)$ and the nearest left-cousin of $u$ in $T^{*}$ is the nearest left-cousin of $u$ in $T$. It follows that \labelcref{T1,T2,T4} hold in $T^{*}$, since they hold in $T$. Now, suppose $u$ is not the right-most vertex in its level in $T^{*}$. Then for any vertex $u'\in V(T^{*})$ such that $L_{T^{*}}(u')\geq L_{T^{*}}(u)+s$, we have $L_{T}(u')\geq L_{T}(u)+s$, which implies $\nu_{T^{*}}(u')=\nu_{T}(u')\geq\nu_{T}(u)=\nu_{T^{*}}(u)$. Thus, \labelcref{T3} holds in $T^{*}$, and therefore $T^{*}$ is a $(K,s)$-tree.
    
    Now suppose that $u_{1}^{*}=u_{t}^{*}$ or $u_{1}^{*}=l(u_{t}^{*})$. It is easy to see that the edges of type \ref{G1}--\ref{G3} in $G^h_{T^{*}}$ map to edges of type \ref{G1}--\ref{G3} in $G^h_T[Q \cup \{u_{t}^*\}]$. Consider an edge $uw$ in $G^h_{T^*}$ added by \ref{G4}. Without loss of generality, we can assume that $L_{T^{*}}(w)\leq L_{T^{*}}(u)$, so $\overrightarrow{uw}\in\overrightarrow{G}_{T^*}^{h}$. Let $u_{a}$ be the $h$-th ancestor of $u$ in $T^{*}$. Let $u_{l}:=l(u_{a})$ and let $u_{r}:=r(u_{a})$ in $T^{*}$.
    
    Suppose $u_{a}=v$. Since $u_{1}^{*}=u_{t}^{*}$ or $u_{1}^{*}=l(u_{t}^{*})$, $D_{T^{*}}(u_{a})\subseteq D_{T}(l(u_{t}^{*}))\cup D_{T}(u_{t}^{*})$. Since $L_T(w)\leq L_T(u)$, it follows that $\overrightarrow{uw}$ is added to $\overrightarrow{G}_{T}^{h}$ by \ref{G4}. If $L_{T^{*}}(u_{a})>0$, then $u_{a},u_{l},u_{r}$ and $u$ are vertices of $T$. In particular, $u_{a}$ is the $h$-th ancestor of $u$ in $T$, $u_{l}$ is the nearest left-cousin of $u_{a}$ in $T$, and $u_{r}$ is the nearest right-cousin of $u_{a}$ in $T$. Since $L_T(w)\leq L_T(u)$, $\overrightarrow{uw}$ is added to $\overrightarrow{G}_{T}^{h}$ by \ref{G4}.
\end{proof}

\subsection{Universality of $G_{T}^{2}$}\label{sec: construction}

We now prove the main theorem of this section. 

\begin{theorem}\label{thm: universal for trees}
    For any $n$-vertex $(4,1)$-tree $T$, the $2$-generated graph $G_{T}^{2}$ is $\mathcal{T}_{n}$-universal.
\end{theorem}

\Cref{thm: universal for trees} is immediate from the following stronger lemma.

\begin{lemma}\label{strong}
    Let $T$ be an $n$-vertex $(4,1)$-tree with root vertex $v$, and let $v_1,\dots, v_{t}$ be the children of $v$ from left to right. Then for any forest $T'$ on $n'\leq n$ vertices, for all $x_1,x_2 \in V(T')$ (not necessarily distinct), there exists a $T$-admissible embedding $\phi: V(T')\to V(G^2_T)$ that satisfies each of the following:
    \begin{enumerate}[label=\upshape \textbf{$\Phi$\arabic{enumi}}, ref=\upshape \textbf{$\Phi$\arabic{enumi}}]
        \item\label{strong1} $\phi(x_1)$ has the minimum level (in $T$) amongst the vertices in $\phi(V(T'))$, and
         \item\label{strong2} 
         if $t=2$ and $n-2\geq n'\geq \nu(v_t)\geq 2$, then $L_T(\phi(x_2))\leq 2$. 
    \end{enumerate}      
\end{lemma}

Let \define{$A_{k}$} be the assertion that \Cref{strong} holds for all trees $T$ with level at most $k$. Let \define{$A_{k,t}$} be the assertion that \Cref{strong} holds for all trees with level at most $k$ whose root has at most $t$ children. We write $(k',t')<(k,t)$ if $k'<k$, or if $k'=k$ and $t'<t$.

\begin{proof}[Proof of \Cref{strong}]

Proceed by induction on $(k,t)$ with the hypothesis that $A_{k,t}$ holds. When $k \leq 2$, $G^2_T$ is isomorphic to a complete graph; it is easy to check that $A_{k}$ holds. Assume now that $k\geq 3$ and that $A_{k',t'}$ holds for all $(k',t')< (k,t)$. We may add edges to $T'$ to make it a tree.



If $\nu(v_t)=1$, then by $A_{k,t-1}$ there is a $T$-admissible embedding $\phi$ of $T'-x_1$ into $G^2_T -v_t$.
If $n'<n$, then extend $\phi$ to a $T$-admissible embedding of $T'$ by mapping $x_1$ to $v_t$. Then \ref{strong1} holds, and \ref{strong2} is vacuous. Now suppose $n'=n$. Let $u':=\phi^{-1}(v)$, and modify $\phi$ by letting $\phi(x_1):=v$ and $\phi(u'):=v_t$. Since $v$ and $v_t$ each see every other vertex in $G^2_T$, $\phi$ is a $T$-admissible embedding satisfying \ref{strong1}. Again, \ref{strong2} is vacuous, implying $\phi$ is a desired embedding. 

We may assume that $\nu(v_t)\geq 2$ for the remainder of the proof. Importantly, \ref{T4} implies that $\nu(v_i)\geq 2$ for all $i\in [t]$. 


If $n'<\nu(v_t)$, then by $A_{k-1}$ there is a $T$-admissible embedding $\phi:V(T')\to D_T[v_t]$ satisfying \ref{strong1}. Since \ref{strong2} is vacuous when $n'<\nu(v_t)$, $\phi$ is a desired embedding of $T'$ into $G_{T}^{2}$. Thus, for the remainder of the proof, we assume that $n'\geq \nu(v_t)$. \newline


\indent {\bf Case 1} ($t=1$){\bf.} Consider the tree $T[D_{T}[v_1]]$ of level $k-1$. If $n'<n$, then by $A_{k-1}$ there is a $T$-admissible embedding $\phi:V(T')\to D_{T}[v_1]$ satisfying \ref{strong1}. Since $t=1$, \ref{strong2} is vacuous, implying $\phi$ is a desired embedding of $T'$ into $G_{T}^{2}$. If $n'=n$, then by $A_{k-1}$ there is a $T$-admissible embedding $\phi:V(T'-x_1)\to D_{T}[v_1]$. Extend $\phi$ to a $T$-admissible embedding of $T'$ into $G_{T}^{2}$ by letting $\phi(x_1):=v$. Then \ref{strong1} holds, and \ref{strong2} is vacuous. \newline

For the remainder of the proof, assume that $t\geq2$. Observe that if $n'\leq \sum_{i=2}^{t} \nu(v_i)$, then by $A_{k,t-1}$ there is a $T$-admissible embedding $\phi$ of $T'$ into $G^2_T-D_{T}[v_{1}]$ satisfying \ref{strong1} and \ref{strong2}. Furthermore, since $n'\geq\nu(v_{t})$, it follows that $\phi$ is a $T$-admissible embedding of $T'$ into $G^2_T$ satisfying \ref{strong1} and \ref{strong2}. Thus, we assume $n' >\sum_{i=2}^{t} \nu(v_i)$ for the remainder of the proof. \newline

\indent {\bf Case 2} ($t=2$){\bf.} Let $u_1,\dots, u_a$ be the children of $v_1$ from left to right, and let $u_{a+1},\dots, u_{a+b}$ be the children of $v_2$ from left to right. Since $\nu(v_1)$ and $\nu(v_2)$ are at least $2$, $a$ and $b$ are positive. We further divide the cases according to the value of $n'$. Since $t=2$, we have $\nu(v_2)\leq n'\leq n$.
\newline

\indent {\bf Case 2.1 }($t=2$ and $\nu(v_2)\leq  n'\leq  n-2$){\bf.} Consider the new tree $T^*$ obtained from $T[D_{T}[u_{1}]\cup\cdots\cup D_{T}[u_{a+b}]]$ by adding a new root $v'$ that has $u_1,\dots, u_{a+b}$ as its children from left to right. By \Cref{lem: tree merge}, $G^2_{T^*}$ is isomorphic to a subgraph of $G^2_T-\{v,v_1\}$. Since $L(T^*)=L(T)-1=k-1$, by $A_{k-1}$ there is a $T^*$-admissible embedding $\phi$ of $T'-x_1$ into $G^2_T-\{v,v_1\}$ such that $\phi(x_2)$ has the minimum level amongst the vertices in $\phi(V(T'-x_1))$. Since $\nu(v_2)\leq  n'\leq  n-2$, $\phi(x_2)$ corresponds to a vertex of level $1$ in $T^{*}$. Thus, $L_{T}(\phi(x_2))=2$. By the minimality of $L_{T}(\phi(x_2))$ in $\phi(V(T'-x_1))$, $v_{2}$ (which corresponds to the root $v'$ of $T^{*}$) is not in $\phi(V(T'-x_1))$. Extend $\phi$ to an embedding of $T'$ into $G_{T}^{2}$ by letting $\phi(x_1):=v_2$. Note that $\phi$ is a $T$-admissible embedding because $n'\geq \nu(v_2)$. Furthermore, since $L_{T}(\phi(x_1))=1$ and $L_{T}(\phi(x_2))=2$, \ref{strong1} and \ref{strong2} hold. Therefore, $\phi$ is a desired embedding of $T'$ into $G_{T}^{2}$. \newline 

\indent {\bf Case 2.2 }($t=2$ and $n'\in \{n-1,n\}$){\bf.} 
Again, we consider the new tree $T^*$ obtained from $T[D_{T}[u_{1}]\cup\cdots\cup D_{T}[u_{a+b}]]$ by adding a new root $v'$ that has $u_1,\dots, u_{a+b}$ as its children from left to right. By \Cref{lem: tree merge}, $G^2_{T^*}$ is isomorphic to a subgraph of $G^2_T-\{v,v_1\}$. Pick a leaf $w$ of $T'-x_{1}$, and let $w'$ be the unique neighbour of $w$ in $T'-x_{1}$. If $T'-x_{1}$ has no leaves, then let $w$ and $w'$ be distinct isolated vertices in $T'-x_{1}$. By \Cref{lem: tree merge}, $G^2_{T^*}$ is isomorphic to a subgraph of $G^2_T-\{v,v_1\}$. Since $L(T^*)=L(T)-1=k-1$, by $A_{k-1}$ there is a $T^*$-admissible embedding $\phi$ of $T'-\{x_1,w\}$ into $G^2_T-\{v,v_1\}$ such that $\phi(w')$ has the minimum level amongst the vertices in $\phi(V(T'-\{x_1,w\}))$. If $n'=n-1$, then $\phi(w')\in\{u_{1},\ldots,u_{a}\}$, in particular, $v_{2}\notin\phi(V(T'-\{x_1,w\}))$. If $n'=n$, then $\phi(w')=v_{2}$. In both cases, $v_{1}$ is adjacent to $\phi(w')$ in $G_{T}^{2}$. Extend $\phi$ to an embedding of $T'$ into $G_{T}^{2}$ by letting
\[\phi(w):= v_1\quad\text{and}\quad\phi(x_1):=
\begin{cases}
    v_2&\text{if } n'=n-1\;,\\
    v&\text{if }n'=n\;.
\end{cases}\]
Then $\phi$ is a $T$-admissible embedding of $T'$ into $G_{T}^{2}$ satisfying \ref{strong1}. Note that \ref{strong2} holds vacuously, so $\phi$ is a desired embedding. \newline

\indent {\bf Case 3 }($t\geq 3$){\bf.} Note that \ref{strong2} is vacuous in this case. Define 
\[\text{$x:= \nu(v_t)$, $y:=\nu(v_{t-1})$ and $z:=\nu(v_{t-2})$}\;.\]
Recall that $x,y,z\geq 2$ and $n'\geq x$. If $x>y$, then by \Cref{lem: feasible critical}, there exists a vertex $w\in V(T')$ and a collection $C_{1},\ldots,C_{p}$ of $w$-components of $T'$ such that $x_{1}\not\in C_{[p]}$ and one of the following cases holds:
\begin{center}
    {\bf Case 3.1.} $C_1,\dots, C_p$ is $(x,y)$-feasible, or$\qquad$ {\bf Case 3.2.} $C_1,\dots, C_p$ is $(x,y)$-critical.
\end{center} 
If $x\leq y$, then an application of \Cref{lem2} to $T'$ with $x-1$ playing the role of $x$ yields a ${(x,y)}$-feasible collection, putting us into Case~3.1 again. Note, it is possible that $w=x_{1}$. For each $i\in [p]$, let $w_i$ be the unique neighbour of $w$ in $C_i$. Let 
\[\text{$T'_0 := T'[ \{w\} \cup C_{[p]}]$ and $T'_1:= T'-V(T'_0)$}\;.\]
Then every edge between $T'_0$ and $T'_1$ is incident to $w$. \newline

We give a brief outline of how to proceed with Case 3.1 and Case 3.2. For Case 3.1, $(x,y)$-feasibility implies certain size conditions on $T'_0$ which allows us to invoke induction to obtain an embedding of $T'_0$ into the graph $2$-generated by the subtree induced by $D_{T}[v_{t-1}]\cup D_{T}[v_t]\cup\{v\}$. By appropriately modifying this embedding we show that \ref{strong1} holds. For Case 3.2, $(x,y)$-criticality implies certain size conditions on $T'_0$ which we use to further find a vertex $w'\in C_1$ and $w'$-components $C'_1,\dots, C'_{p'}$ of $T'_0[C_1]$ with appropriate sizes. This divides $C_1$ into two parts, $C'_{[p']}$ and $C_1\setminus C'_{[p']}$. Using these parts, we define a sequence of pairwise-disjoint forests $T_{p}'',T_{p-1}'',\ldots,T''_{1},T_{0}''$ in $T'$ which we embed one by one (in the given order) into $G^{2}_{T}$ by induction. Again, by induction we embed the remaining part of $T'$ into $G^{2}_{T}$. These embeddings are chosen so that the concatenation $\phi$ is a $T$-admissible embedding of $T'$. In order to guarantee that $\phi$ is an embedding of $T'$, we make careful choices when assigning the image of the vertices $w$, $w'$ and $w_{1}$; there are three subcases to consider. In Case I, $C'_{[p']}$ is assigned to a tree $T''_0$ together with $w$ (see \Cref{fig:Case3-2-I}). In Case II and III, $C_1\setminus C'_{[p']}$ is assigned to a tree $T''_0$ (see \Cref{fig:Case3-2-II,fig:Case3-2-III}). Depending on whether or not $w=x_1$, we carefully choose which vertices of $w,w',w'_1$ belong to $T''_0$; this distinguishes Case II and III.
When we apply induction to each $T''_i$ in the sequence, we make use of \ref{strong2} to control the image of some special vertices (such as $w$ and $w'$). Thus, when embedding each $T''_i$, we keep track of the number of vertices remaining after embedding each piece into $D_{T}[v_t]$ to make sure the conditions to invoke \ref{strong2} are satisfied. \newline

\indent {\bf Case 3.1 } ($C_1,\dots, C_p$ is $(x,y)$-feasible){\bf.}
Note that \(x\leq |T'_0|\leq x+y-2\). Let $Q:=D_T[v_{t-1}]\cup D_T[v_t]$ and consider the tree $T^{*}:=T[Q\cup\{v\}]$. By \Cref{lem: tree merge}, $G^2_{T^*}$ is a subgraph of $G^2_T$. Since $|T'_0|\leq x+y-2$, by $A_{k,2}$ there is a $T^*$-admissible embedding $\phi_0$ of $T'_0$ into $G^2_{T^*}$ such that $\phi_0(w)$ has the minimum level amongst the vertices in $\phi_0(V(T'_0))$. Since $|T'_0|\leq x+y-2<|Q|+1=|T^{*}|$, $v\not\in\phi_0(V(T'_0))$. So in fact, $\phi_0$ is a $T$-admissible embedding of $T'_0$ into $G^2_T$. Moreover, since $x\leq|T'_0|$, by the minimality of $L_{T}(\phi_0(w))$ in $\phi_0(V(T'_0))$, we have $\phi_0(w)=v_t$.

Let $T_{1}:=T-\phi_0(V(T'_0))$. Since $\phi_0(V(T'_0))$ is $T$-admissible, and $x\leq|T'_0|\leq x+y-2$, $T_{1}$ is a $(4,1)$-tree in which the root $v$ has exactly $t-1$ children. By $A_{k,t-1}$, there is a $T_1$-admissible embedding $\phi_1$ of $T'_1$ into $G^2_T-\phi_0(V(T'_0))$ such that if $w\not=x_{1}$, then $\phi_{1}(x_{1})$ has the minimum level amongst the vertices in $\phi_{1}(V(T_{1}))$. Since $\phi_0(w)=v_{t}$ sees every other vertex in $G^2_T$, combining $\phi_{0}$ and $\phi_{1}$ yields a $T$-admissible embedding $\phi$ of $T'$ into $G^2_T$. It remains to check that \ref{strong1} holds.

If $w=x_{1}$, then the only way \ref{strong1} does not hold is when $n'=n$. In this case, we modify $\phi$ by mapping the preimage of $v$ to $v_{t}$, and mapping $x_{1}$ to $v$. Since $v$ and $v_{t}$ each see every other vertex of $G^2_T$, $\phi$ is a desired embedding satisfying \ref{strong1}. 

If $w\not=x_{1}$, then $\phi_{1}(x_{1})$ has the minimum level amongst the vertices in $\phi_{1}(V(T_{1}))$. Since we are assuming $n' >\sum_{i=2}^{t} \nu(v_i)$, it follows that $L_{T}(\phi_{1}(x_{1}))\leq1$, implying \ref{strong1} holds. \newline

\indent{\bf Case 3.2 } ($C_1,\dots, C_p$ is $(x,y)$-critical){\bf.} 
In this case, we further partition $T'_0$ into several forests. We may assume that $y<x$ and that $|C_1|\leq\cdots \leq |C_p|$. Consider the following claim.
\begin{claim}\label{claim: T'0-C' sizes}
There exists a vertex $w'\in C_1$ and a collection $C'_1,\dots, C'_{p'}$ of $w'$-components of $T'[C_1]$ such that $w_1\notin C'_{[p']}$ and
\begin{enumerate}[label=(\alph*)]
    \item\label{case I} $x\leq|T'_0-(C'_{[p']}\cup\{w\})|\leq x+y-2$, or
    \item\label{case II} $x\leq|C_{[2,p]}\cup\{w, w'\}\cup C'_{[p']}|\leq x+y-2$.
\end{enumerate}
\end{claim}
\begin{proof}
Let $z':=|T'_0|-x-y$. By \ref{T1}, $y+z\geq x$. Furthermore, since $C_1,\dots, C_p$ is $(x,y)$-critical, $x+y-1\leq |T'_0|\leq 2x-2$. Therefore, $-1\leq z'\leq z-2$. We proceed by cases on the value of $z'$.

First, consider the case when $z'\leq y-2$. If $z'+1=0$, then choosing any vertex $w'\in C_1$ and the empty collection of $w'$-components shows that \labelcref{case I} holds. If $z'+1\geq1$, then since $|C_{1}|\geq y> z'+1$, applying \Cref{lem2} with $z'+1$ playing the role of $x$ yields a vertex $w'\in C_1$ and collection of $w'$-components $C'_1,\dots, C'_{p'}$ of $T'[C_1]$ with $w_1\notin C'_{[p']}$ that satisfies:
\[z'+1\leq|C'_{[p']}|\leq2(z'+1)-1\leq y+z'-1\;.\]
Since $|T'_0|=x+y+z'$, \labelcref{case I} holds.

Now, consider the case when $z'\geq  y-1$. Since $p\geq 2$ and $|C_1|\leq |C_2|$, 
$$2|C_{[2,p]}|+1 \geq |C_1\cup C_{[2,p]}\cup \{w\}|= |T'_0|=x+y+z'\;.$$ 
Furthermore, since $T$ is a $(4,1)$-tree, by \ref{T2}, $x\leq 4y$. Since $z'\geq y-1$, $z'\geq x-3y-1$. Therefore,
$|C_{[2,p]}|\geq  \frac{x+y+z'-1}{2} \geq x-y-1$, implying
\[2(x-|C_{[2,p]}|-2)\leq x+y-|C_{[2,p]}|-3\;.\]
Moreover, since $|C_1|\geq y>x-|C_{[2,p]}|-2$, applying \Cref{lem2} with $x-|C_{[2,p]}|-2$ playing the role of $x$ yields a vertex $w'\in C_1$ and a collection $C'_{1}, \dots, C'_{p'}$ of $w'$-components of $T'[C_{1}]$ with $w_1\notin C'_{[p']}$ that satisfies:
\[x-|C_{[2,p]}|-2\leq|C'_{[p']}|\leq2(x-|C_{[2,p]}|-2)-1\leq  x+y-|C_{[2,p]}|-4\;.\] 
Since $|C_{[2,p]}\cup\{w, w'\}\cup C'_{[p']}|=|C_{[2,p]}|+2+|C'_{[p']}|$, the extremes of the above inequality implies that \labelcref{case II} holds. This proves the claim.
\end{proof}

Let $w'$ and $C'_1,\dots, C'_{p'}$ be as in the above claim. Then either \labelcref{case I} or \labelcref{case II} holds. Let $C':=C'_{[p']}$ and define
\[C_0:=
\begin{cases}
    \{w\}&\text{ if }|T'|\geq x+y+z\;,\\
    V(T')\setminus C_{[p]}&\text{ if } |T'|\leq x+y+z-1\;.
\end{cases}\]
Note that $w\in C_{0}$. We now define a sequence of pairwise-disjoint forests $T''_{0},T_{1}'',\ldots,T_{p}''$ depending on the following cases.\newline 

\indent {\bf Case I.} If \labelcref{case I} holds, then let 
\[\text{$T''_0:= T'[C'\cup C_0]$,  $T''_1:= T'_0[C_1]-C' $ and $T''_j := T'[C_j]$ for $j\geq 2$}\;.\]

\indent {\bf Case II.} If \labelcref{case II} holds and $w\neq x_1$, then let 
\[\text{$T''_0:= T'[C_1\cup C_0]-(C'\cup\{w, w'\})$, $T''_1:= T'_0[C'\cup\{w, w'\}]$ and $T''_j := T'[C_j]$ for $j\geq 2$}\;.\]

\indent {\bf Case III.} If \labelcref{case II} holds and $w=x_1$, then pick a leaf $z$ of $T'[C_1\cup C_0]-(C'\cup\{w, w'\})$, and let $z'$ be the unique neighbour of $z$. If $T'[C_1\cup C_0]-(C'\cup\{w, w'\})$ has no leaves, then let $z$ and $z'$ be distinct isolated vertices in $T'[C_1\cup C_0]-(C'\cup\{w, w'\})$. Let
\[\text{$T''_0:= T'[C_1\cup C_0]-(C'\cup\{z, w'\})$, $T''_1:= T'_0[C'\cup\{z, w'\}]$ and $T''_j := T'[C_j]$ for $j\geq 2$}\;.\]

\begin{figure}[!ht]
    \centering

    \begin{tikzpicture}[/tikz/xscale=1, /tikz/yscale=1]

    \begin{scope}[every node/.style={ellipse, very thick,draw, lightgray}]
        \node[minimum width=6cm, minimum height=5cm] (C1) at (0,0) {};
        \node[minimum width=1.5cm, minimum height=5cm] (C2) at (4.5,0) {\textcolor{black}{$T_{2}''$}};
        \node[minimum width=1.5cm, minimum height=5cm] (Cp) at (4.5+3,0) {\textcolor{black}{$T_{p}''$}};

        \node[minimum width=4cm, minimum height=2cm] (Top) at (0,1.2) {};
        \node[minimum width=1cm, minimum height=1.5cm] (Top1) at (-1,1.2) {};
        \node[minimum width=1cm, minimum height=1.5cm] (Top2) at (1,1.2) {};
        
        \node[minimum width=4cm, minimum height=2cm] (Bot) at (0,-1.2) {};
        \node[minimum width=1cm, minimum height=1.5cm] (Bot1) at (-1,-1.2) {};
        \node[minimum width=1cm, minimum height=1.5cm] (Bot2) at (1,-1.2) {};
    \end{scope}

    \begin{scope}[every node/.style={draw, shape = circle, fill = black, minimum size = 0.1cm, inner sep=1.5pt}]
        \node(w) at (1.9,3.1) {};
        \node(w') at (0,0.5) {};
        \node(w1) at (-1,1.7) {};
        \node(w2) at (4.5,2) {};
        \node(wp) at (4.5+3,2) {};

        \node(v1) at (-1,0.7) {};
        \node(v2) at (1,0.7) {};
        \node(v3) at (-1,-0.7) {};
        \node(v4) at (1,-0.7) {};
    \end{scope}

    \begin{scope}[every edge/.style={draw, thick}]
        \path [-] (w) edge (w1);
        \path [-] (w) edge (w2);
        \path [-] (w) edge (wp);

        \path [-] (w') edge (v1);
        \path [-] (w') edge (v2);
        \path [-] (w') edge (v3);
        \path [-] (w') edge (v4);
    \end{scope}

    \begin{scope}
        \node (C1) at (-2.2,-2.2) {$C_{1}$};
        \node (C2) at (4.5+0.8,-2.2) {$C_{2}$};
        \node (Cp) at (4.5+3+0.8,-2.2) {$C_{p}$};

        \node (TopDots) at (0,1.5) {$\cdots$};
        \node (BotDots) at (0,-1.5) {$\cdots$};
        \node (MidDots) at (4.5+1.5,0) {$\cdots$};

        \node[shape = circle, fill=white, minimum size = 0.1cm, inner sep=1.5pt] (T) at (1.9,0.8) {$T_{1}''$};
        \node[shape = circle, fill=white, minimum size = 0.1cm, inner sep=1.5pt] (C) at (1.9,-0.8) {$C'$};
        \node (C1') at (-1,-1.5) {$C_{1}'$};
        \node (Cp') at (1,-1.5) {$C_{p'}'$};

        \node (w) at (1.9,3+0.36) {$w$};
        \node (w') at (0,0.5+0.35) {$w'$};
        \node (w1) at (-1,1.7-0.35) {$w_{1}$};
        \node (w2) at (4.5,2-0.35) {$w_{2}$};
        \node (wp) at (4.5+3,2-0.35) {$w_{p}$};
    \end{scope}
    \end{tikzpicture}
    
    \caption{Structures of $T'_{0}$ in Case I.}
    \label{fig:Case3-2-I}
\end{figure}
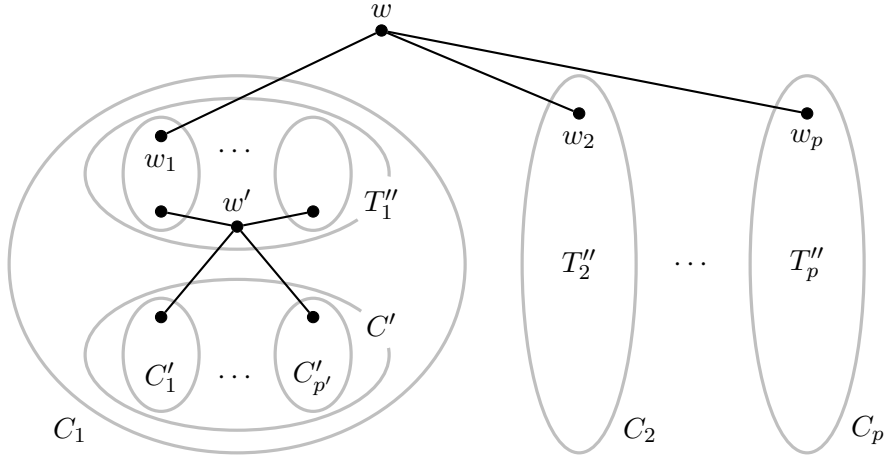

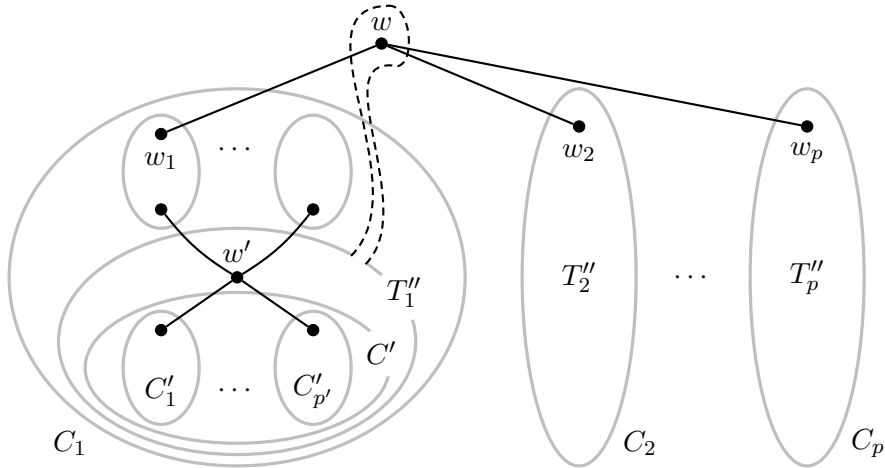
\begin{figure}[!ht]
    \centering

    \begin{tikzpicture}[/tikz/xscale=1, /tikz/yscale=1]

    \begin{scope}[every node/.style={ellipse, very thick,draw, lightgray}]
        \node[minimum width=6cm, minimum height=5cm] (C1) at (0,0) {};
        \node[minimum width=1.5cm, minimum height=5cm] (C2) at (4.5,0) {\textcolor{black}{$T_{2}''$}};
        \node[minimum width=1.5cm, minimum height=5cm] (Cp) at (4.5+3,0) {\textcolor{black}{$T_{p}''$}};
        
        \node[minimum width=1cm, minimum height=1.5cm] (Top1) at (-1,1.4) {};
        \node[minimum width=1cm, minimum height=1.5cm] (Top2) at (1,1.4) {};

        \node[minimum width=4.7cm, minimum height=3cm] (BotOuter) at (0,-0.85) {};
        \node[minimum width=4cm, minimum height=2cm] (Bot) at (0,-1.2) {};
        \node[minimum width=1cm, minimum height=1.5cm] (Bot1) at (-1,-1.2) {};
        \node[minimum width=1cm, minimum height=1.5cm] (Bot2) at (1,-1.2) {};
    \end{scope}

    \begin{scope}[every node/.style={draw, shape = circle, fill = black, minimum size = 0.1cm, inner sep=1.5pt}]
        \node(w) at (1.9,3.1) {};
        \node(w') at (0,0) {};
        \node(w1) at (-1,1.9) {};
        \node(w2) at (4.5,2) {};
        \node(wp) at (4.5+3,2) {};

        \node(v1) at (-1,0.9) {};
        \node(v2) at (1,0.9) {};
        \node(v3) at (-1,-0.7) {};
        \node(v4) at (1,-0.7) {};
    \end{scope}

    \begin{scope}
        \node[shape = circle, fill=white, minimum size = 0.21cm, inner sep=1.5pt] (patch) at (1.607,0.24) {};
    \end{scope}

    \begin{scope}[every edge/.style={draw, thick}]
        \path [-] (w) edge (w1);
        \path [-] (w) edge (w2);
        \path [-] (w) edge (wp);

        \path [-] (w') edge[bend left=10] (v1);
        \path [-] (w') edge[bend right=10] (v2);
        \path [-] (w') edge (v3);
        \path [-] (w') edge (v4);

        \path[densely dashed, thick, draw, black] (1.5,0.285) edge [in=180] (1.9+0.1,3.1+0.5);
        \path[densely dashed, thick, draw, black] (1.5+0.2-0.01,0.285-0.11) edge [in=170] (1.9+0.22,3.1-0.3);
        \path[densely dashed, thick, draw, black, out=5, in=-10, looseness=1] (1.9+0.1-0.01,3.1+0.5) to (1.9+0.2-0.04,3.1-0.3);
    \end{scope}

    \begin{scope}
        \node (C1) at (-2.2,-2.2) {$C_{1}$};
        \node (C2) at (4.5+0.8,-2.2) {$C_{2}$};
        \node (Cp) at (4.5+3+0.8,-2.2) {$C_{p}$};

        \node (TopDots) at (0,1.7) {$\cdots$};
        \node (BotDots) at (0,-1.5) {$\cdots$};
        \node (MidDots) at (4.5+1.5,0) {$\cdots$};

        \node[shape = circle, fill=white, minimum size = 0.1cm, inner sep=1.5pt] (T) at (2.2,-0.2) {$T_{1}''$};
        \node[shape = circle, fill=white, minimum size = 0.1cm, inner sep=1.5pt] (C) at (1.9,-1) {$C'$};
        \node (C1') at (-1,-1.5) {$C_{1}'$};
        \node (Cp') at (1,-1.5) {$C_{p'}'$};

        \node (w) at (1.9,3+0.36) {$w$};
        \node (w') at (0,0.35) {$w'$};
        \node (w1) at (-1,1.7-0.15) {$w_{1}$};
        \node (w2) at (4.5,2-0.35) {$w_{2}$};
        \node (wp) at (4.5+3,2-0.35) {$w_{p}$};

    \end{scope}
    \end{tikzpicture}
    
    \caption{Structures of $T'_{0}$ in Case II.}
    \label{fig:Case3-2-II}
\end{figure}
    
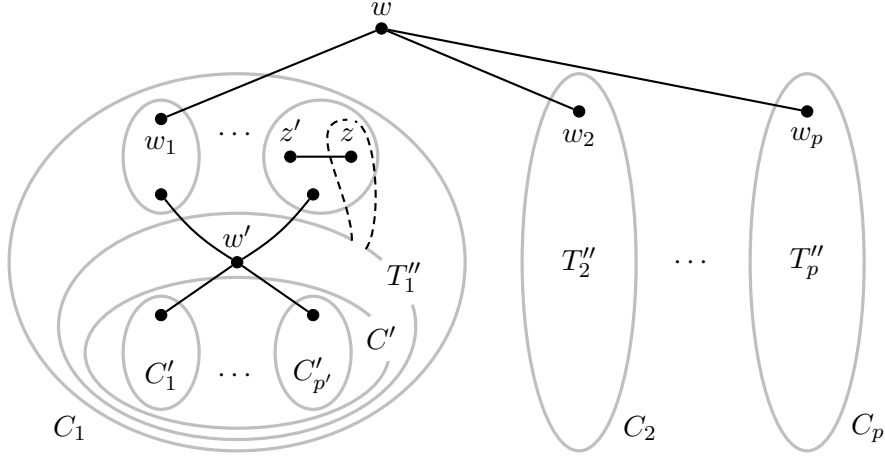
\begin{figure}[!ht]
    \centering

    \begin{tikzpicture}[/tikz/xscale=1, /tikz/yscale=1]

    \begin{scope}[every node/.style={ellipse, very thick,draw, lightgray}]
        \node[minimum width=6cm, minimum height=5cm] (C1) at (0,0) {};
        \node[minimum width=1.5cm, minimum height=5cm] (C2) at (4.5,0) {\textcolor{black}{$T_{2}''$}};
        \node[minimum width=1.5cm, minimum height=5cm] (Cp) at (4.5+3,0) {\textcolor{black}{$T_{p}''$}};
        
        \node[minimum width=1cm, minimum height=1.5cm] (Top1) at (-1,1.4) {};
        \node[minimum width=1.5cm, minimum height=1.5cm] (Top2) at (1.1,1.4) {};

        \node[minimum width=4.7cm, minimum height=3cm] (BotOuter) at (0,-0.85) {};
        \node[minimum width=4cm, minimum height=2cm] (Bot) at (0,-1.2) {};
        \node[minimum width=1cm, minimum height=1.5cm] (Bot1) at (-1,-1.2) {};
        \node[minimum width=1cm, minimum height=1.5cm] (Bot2) at (1,-1.2) {};
    \end{scope}

    \begin{scope}[every node/.style={draw, shape = circle, fill = black, minimum size = 0.1cm, inner sep=1.5pt}]
        \node(w) at (1.9,3.1) {};
        \node(w') at (0,0) {};
        \node(w1) at (-1,1.9) {};
        \node(w2) at (4.5,2) {};
        \node(wp) at (4.5+3,2) {};

        \node(v1) at (-1,0.9) {};
        \node(v2) at (1,0.9) {};
        \node(v3) at (-1,-0.7) {};
        \node(v4) at (1,-0.7) {};

        \node(z') at (1.1-0.4,1.4) {};
        \node(z) at (1.1+0.4,1.4) {};
    \end{scope}

    \begin{scope}
        \node[shape = circle, fill=white, minimum size = 0.21cm, inner sep=1.5pt] (patch) at (1.607,0.24) {};
    \end{scope}

    \begin{scope}[every edge/.style={draw, thick}]
        \path [-] (w) edge (w1);
        \path [-] (w) edge (w2);
        \path [-] (w) edge (wp);

        \path [-] (w') edge[bend left=10] (v1);
        \path [-] (w') edge[bend right=10] (v2);
        \path [-] (w') edge (v3);
        \path [-] (w') edge (v4);
        \path [-] (z) edge (z');

        
        \path[densely dashed, thick, draw, black, looseness=1] (1.51,0.285) edge [in=180, out=90] (1.1+0.35,1.4+0.49);
        \path[densely dashed, thick, draw, black, out=50, in=-20, looseness=0.5] (1.5+0.2-0.01,0.285-0.11) to (1.1+0.45,1.4+0.47);
    \end{scope}

    \begin{scope}
        \node (C1) at (-2.2,-2.2) {$C_{1}$};
        \node (C2) at (4.5+0.8,-2.2) {$C_{2}$};
        \node (Cp) at (4.5+3+0.8,-2.2) {$C_{p}$};

        \node (TopDots) at (0,1.7) {$\cdots$};
        \node (BotDots) at (0,-1.5) {$\cdots$};
        \node (MidDots) at (4.5+1.5,0) {$\cdots$};

        \node[shape = circle, fill=white, minimum size = 0.1cm, inner sep=1.5pt] (T) at (2.2,-0.2) {$T_{1}''$};
        \node[shape = circle, fill=white, minimum size = 0.1cm, inner sep=1.5pt] (C) at (1.9,-1) {$C'$};
        \node (C1') at (-1,-1.5) {$C_{1}'$};
        \node (Cp') at (1,-1.5) {$C_{p'}'$};

        \node (w) at (1.9,3+0.36) {$w$};
        \node (w') at (0,0.35) {$w'$};
        \node (w1) at (-1,1.7-0.15) {$w_{1}$};
        \node (w2) at (4.5,2-0.35) {$w_{2}$};
        \node (wp) at (4.5+3,2-0.35) {$w_{p}$};

        \node(z') at (1.1-0.4,1.4+0.35) {$z'$};
        \node(z) at (1.1+0.35,1.4+0.29) {$z$};
    \end{scope}
    \end{tikzpicture}
    
    \caption{Structures of $T'_{0}$ in Case III.}
    \label{fig:Case3-2-III}
\end{figure}

Then in each of the cases, $\{V(T_{0}''),V(T_{1}''),\ldots,V(T_{p}''),V(T_{1}')\setminus C_{0}\}$ is a partition of $V(T')$. Note that $T_{1}'-C_{0}=T_{1}'$ when $|T'|\geq x+y+z$, and $T_{1}'-C_{0}$ is empty when $|T'|\leq x+y+z-1$. Furthermore, $w',w_1\in T''_1$ in Case I, $w\in T''_1$ in Case II, and either $w_1=w'$ or $w_1\notin T''_1$ in Case III. Moreover, \Cref{claim: T'0-C' sizes} implies that in each of the three cases we have:
\begin{align}\label{eq: T'' sizes}
 x\leq |T''_1|+|T''_{[2, p]}| \leq x+y-2\;.
\end{align}
We proceed by embedding $T''_p, \dots, T''_2, T''_1, T''_0$ and $T'_1-C_0$ in the given order. 
First, we embed $T''_p, \dots, T''_2$ recursively. Let $T^p:=T$ and let $\phi_j:V(T''_j)\to V(T^j)$ be a $T^j$-admissible embedding such that $\phi_{j}(w_j)$ has the minimum level amongst the vertices in $\phi_j(V(T''_j))$, where $T^j=T-\bigcup_{j'=j+1}^{p}\phi_{j'}(V(T''_{j'}))$. We explain how to obtain such a sequence of embeddings.

 By $(x, y)$-criticality, $|C_{[2, p]}|\leq x-2$, which implies that the image of a $T$-admissible embedding $C_{[2, p]}\to V(G^{2}_{T})$ is contained entirely in $D_T[v_t]$. Suppose that for a given $j\geq2$, $T''_{j'}$ is already embedded for all $j'>j$. Let $u'_1,\dots, u'_{a}$ be the children of $v_t$ in $T^j$ and consider the subtree induced by $D_{T^j}[v_t]$. Since $T$ is a $(4,1)$-tree and $\nu_{T}(v_{t-1})=y$, \ref{T3} implies that $\nu_{T^j}(u'_j)\leq\nu_{T}(u'_j)\leq y$ for all $i\in [a]$. Furthermore, by $(x,y)$-criticality, $|T''_j|\geq y$. Hence, by the induction hypothesis $A_{k-1}$, we obtain a $T^j$-admissible embedding $\phi_{j}:V(T''_{j})\to D_{T^j}[v_t]$ such that $\phi_{j}(w_j)$ has the minimum level amongst the vertices in $\phi_{j}(V(T''_{j}))$. Moreover, since $|T''_j|\geq \nu_{T^j}(u'_a)$, $\phi_{j}(w_j)$ is a child of $v_t$. \newline

Second, we embed $T''_1$. 
Let $T^*$ be the tree obtained from $T[D_T[v_{t-1}]\cup D_T[v_{t}]] -\bigcup_{i=2}^{p}\phi_i(V(T''_{i}))$ by adding a new root $v'$. By the definition of $(x,y)$-critical collection, $|C_{[2, p]}| \leq  x-2 = \nu_T(x_t)-2$, which implies $\nu_{T^*}(v_t)\geq 2$. Let $Q_1:=(D_T[v_{t-1}]\cup D_T[v_{t}])\setminus\bigcup_{i=2}^{p}\phi_i(V(T''_{i}))$. Then  \Cref{lem: tree merge} says that  $G^{2}_{T^*}$ is isomorphic to a subgraph of $G^{2}_T[Q_1\cup \{v\}]$. Furthermore, \eqref{eq: T'' sizes} implies 
 \[\nu_{T^*}(v_{t}) = x - |T''_{[2, p]}| \leq  |T''_{[p]}|-|T''_{[2, p]}|= |T''_1|\leq x+y-2 - |T''_{[2, p]}| = \nu_{T^*}(v_t)+ \nu_{T^*}(v_{t-1})-2\;.\]
Thus, $|T^*|-2\geq |T''_1|\geq \nu_{T^*}(v_t)\geq 2$. We apply the induction hypothesis $A_{k,2}$ to obtain a $T$-admissible embedding $\phi_1$ of $T''_1$ into $G^2_T-\bigcup_{i=2}^{p}\phi_i(V(T''_{i}))$ where $w'$ plays the role of $x_1$, and depending on the cases, the following vertex plays the role of $x_2$.
\begin{center}
\begin{tabular}{ |c|c|c|c|c|c|c|c|c|c| }
\hline
   &  Case I & Case II & Case III\\ \hline 
    Role of $x_2$ & $w_1$ & $w$ & $z$ \\
\hline
\end{tabular}
\end{center}
Therefore, $\phi_1(w')=v_{t}$, and depending on the cases, $\phi_1(w_1)$ or $\phi_1(w)$ or $\phi_1(z)$ has level most $2$. \newline 
 
Third, we embed $T''_0$. Let $T^{**}$ be the tree obtained from $T[D_T[v_{t-2}]\cup D_{T}[v_{t-1}]]- \bigcup_{i=1}^{p}\phi_i(V(T''_{i}))$ by adding a new root $v'$. Since $|T''_1| \leq \nu_{T^*}(v_t)+\nu_{T^*}(v_{t-1})-2$, $\nu_{T^{**}}(v_{t-1})\geq 2$.
Let $Q_0:=(D_T[v_{t-2}]\cup D_T[v_{t-1}])\setminus\bigcup_{i=1}^{p}\phi_i(T''_{i})$. Then \Cref{lem: tree merge} says that $G^{2}_{T^{**}}$ is isomorphic to a subgraph of $G^{2}_T[Q_0\cup\{v\}]$.  
Furthermore, by the definition of $C_0$ and by $(x,y)$-criticality,
\[|T''_0|+|T''_1|+|T''_{[2, p]}|\leq |T'_0|+|C_0|-1 \leq \max\{ 2x-2, x+y+z-1\} = x+y+z-1\]
where the last equality holds since $x\leq y+z$.
Thus, $|T^{**}|-2\geq |T''_0| \geq \nu_{T^{**}}(v_{t-1})\geq 2$. We apply the induction hypothesis $A_{k,2}$ to embed $T''_0$ into $G^2_T[Q_0\cup\{v\}]$ by considering the following cases.

\begin{enumerate}
    \item[\scriptsize$\bullet$] If $x_1\in C_0\setminus\{w\}$, then Case III does not occur, since $x_1\not= w$. We apply the induction hypothesis $A_{k,2}$
    to obtain a $T$-admissible embedding $\phi_0$ of $T''_0$ into $G^2_T[Q_0\cup\{v\}]$ with the following vertices playing the roles of $x_1$ and $x_2$ respectively.
    \begin{center}
    \begin{tabular}{ |c|c|c|c|c|c|c|c|c|c| }
    \hline
       &  Case I & Case II \\ \hline 
        Role of $x_1$ & $x_1$ & $x_1$  \\ \hline
        Role of $x_2$ & $w$ & $w_1$  \\
    \hline
    \end{tabular}
    \end{center}
    Therefore, $\phi_0(x_1)=v_{t-1}$, and depending on the cases, $\phi_0(w)$ or $\phi_0(w_1)$ has level at most $2$.

    \item[\scriptsize$\bullet$] If $x_1\notin C_0\setminus\{w\}$, then either $x_{1}\in V(T'_1)\setminus C_0$ or $w=x_1$. In both cases, we apply the induction hypothesis $A_{k,2}$ to obtain a $T$-admissible embedding $\phi_0$ of $T''_0$ into $G^2_T[Q_0\cup\{v\}]$ with the following vertices playing the roles of $x_1$ and $x_2$ respectively.
    \begin{center}
    \begin{tabular}{ |c|c|c|c|c|c|c|c|c|c| }
    \hline
       &  Case I & Case II & Case III \\ \hline 
        Role of $x_1$ & $w$ & $w_1$ & $w$  \\ \hline
        Role of $x_2$ &  &  & $z'$ \\
    \hline
    \end{tabular}
    \end{center}
    Therefore, $\phi_0(w)=v_{t-1}$ or $\phi_0(w_1)=v_{t-1}$, and in Case III, the level of $\phi_0(z')$ is at most $2$.
\end{enumerate}

Finally, we embed $T'_1-C_0$. If $T'_1-C_0$ is empty, let $\phi'$ be the empty embedding. 
Otherwise, consider the tree $T^{\#}:=T - \bigcup_{i=0}^{p} \phi_{i}(V(T''_i))$. By the induction hypothesis $A_{k,t-2}$, there is a $T$-admissible embedding $\phi'$ of $T'_1-C_0$ into $G_{T^{\#}}^2 = G^2_T - \bigcup_{i=0}^{p} \phi_{i}(V(T''_i))$  such that $\phi'(x_1)$ has the minimum level amongst the vertices in $\phi'(V(T'_1)\setminus C_0)$ when $x_1\in V(T'_1)\setminus C_0$. Note that when $T'_1- C_0$ is not empty, $|T'|\geq x+y+z$, which implies that $v_{t-2}\in\phi'(V(T'_1)\setminus C_0)$. So $\phi'(x_1)$ has level at most $1$.

\begin{figure}[!h]
    \centering
    \hspace*{2.3cm}
    \begin{tikzpicture}
    \begin{scope}
        \draw[very thick, draw=lightgray] (-1.2-3,-0.5) -- (1.2-3,-0.5) -- (-3,2) -- (-1.2-3,-0.5);

        \draw[very thick, draw=lightgray] (-1.2,-0.5) -- (1.2,-0.5) -- (0,2) -- (-1.2,-0.5);

        \draw[very thick, draw=lightgray] (-1.2+3,-0.5) -- (1.2+3,-0.5) -- (3,2) -- (-1.2+3,-0.5);
    \end{scope}
    
    \begin{scope}[draw, shape = circle, fill = black, minimum size = 0.1cm, inner sep=1.6pt]
        \node[draw, shape = circle, fill = black, minimum size = 0.1cm, inner sep=1.6pt] (v) at (0,3.5) {};
        \node[draw, shape = circle, fill = black, minimum size = 0.1cm, inner sep=1.6pt] (vleft) at (-3,2) {};
        \node[draw, shape = circle, fill = black, minimum size = 0.1cm, inner sep=1.6pt] (vmiddle) at (0,2) {};
        \node[draw, shape = circle, fill = black, minimum size = 0.1cm, inner sep=1.6pt] (vright) at (3,2) {};
        
        \node (node) at (0,3.9) {$v$};
        \node (node left) at (-3-0.7,2) {$v_{t-2}$};
        \node (node middle) at (0-0.7,2) {$v_{t-1}$};
        \node (node right) at (3-0.7,2) {$v_{t}$};

        \node (text left) at (-1.5,-1.25) {$\phi_{0}(V(T_{0}''))$};
        \node (text middle) at (1.5,-1.25) {$\phi_{1}(V(T_{1}''))$};
        \node (text right) at (5.35,1) {$\bigcup_{i=2}^{p}\phi_{i}(V(T_{i}''))$};
    \end{scope}
    
    \begin{scope}[every edge/.style={draw, thick}]
        \path [-] (v) edge (vleft);
        \path [-] (v) edge (vmiddle);
        \path [-] (v) edge (vright);
    \end{scope}

    \begin{scope}[every node/.style={ellipse, thick, draw, black}]
        \node[densely dashed, minimum width=4cm, minimum height=2.5cm, rotate=60] (left) at (-1,1) {};
        \node[densely dashed, minimum width=4cm, minimum height=2.5cm, rotate=60] (middle) at (2,1) {};
        \node[densely dashed, minimum width=1.5cm, minimum height=1.5cm] (right) at (4,0-0.1) {};
    \end{scope}
    \end{tikzpicture}\vspace{-2pc}
    \caption{Each tree $T''_i$ embeds as above. In addition, $v_t= \phi_1(w')$ and $v_{t-1}\in \{\phi_0(x_1), \phi_0(w)\}$ in Cases I and III, and $v_{t-1}\in \{\phi_0(x_1), \phi_0(w_1)\}$ in Case II.}
\end{figure}
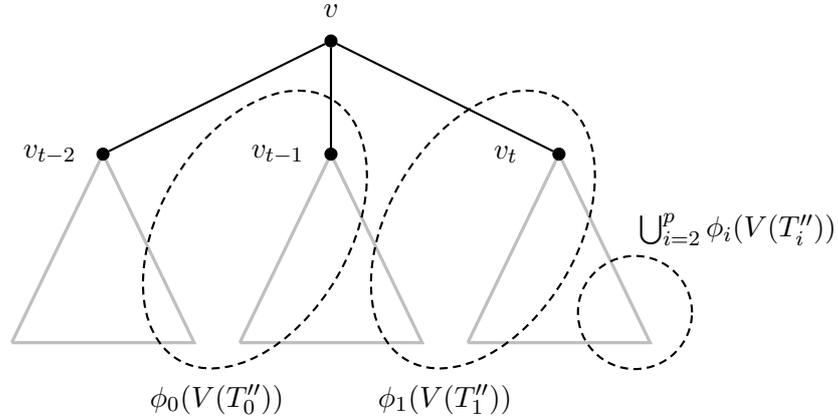

The final embedding $\phi:V(T')\to V(G^{2}_{T})$ is obtained by concatenating $\phi_p, \dots, \phi_1, \phi_0$ and $\phi'$; we check that $\phi$ is indeed an embedding.
Note that any edge between $T''_0, T''_1, \dots, T''_p$ and $T'_1-C_0$ is incident to one of the three vertices $w, w'$ and $z$ (where vertex $z$ is only considered in Case III). Since $\phi(w')=v_{t}$ sees every other vertex in $G^2_T$, the edges incident with $w'$ are mapped to edges in $G^2_T$. 

We check that $\phi(w)$ and $\phi(w_i)$ are adjacent for each $i$. In every case, each of the vertices $\phi(w),\phi(w_2),\ldots,\phi(w_p)$ have level at most $2$, implying they are all pairwise-adjacent by \Cref{newdef}.
The only case in which the level of $\phi(w_1)$ is larger than $2$ is Case III. In Case III, $w$ and $w_1$ belong to $T''_0$, implying they are adjacent since $\phi_0$ is an embedding. 
Now consider the edges incident with $z$. In Case III, $zw, zw', zz'$ are the only such edges possible in $T'$. Since each of the vertices $\phi(w),\phi(w'),\phi(z)$ and $\phi(z')$ have level at most $2$, they are all pairwise-adjacent by \Cref{newdef}. This shows that $\phi$ is indeed a $T$-admissible embedding.

It remains to check that $\phi$ satisfies \ref{strong1}. Note, either $x_1\in V(T'_1)\setminus C_0$ or $x_1\in C_0\setminus\{w\}$ or $x_1=w$. We proceed by cases.

\begin{enumerate}
    \item[\scriptsize$\bullet$] If $x_1\in V(T'_1)\setminus C_0$, \ref{strong1} holds by the choice of $\phi'$.
    \item[\scriptsize$\bullet$]  If $x_1\in C_0\setminus\{w\}$, then $\phi(x_1)=v_{t-1}$ has level $1$. Since $C_0\setminus\{w\}\not=\varnothing$, $|T'|\leq x+y+z-1$, which implies $v\notin \phi(V(T'))$. Therefore, \ref{strong1} holds. 
    \item[\scriptsize$\bullet$] If $x_1=w$ and $n'<n$, then either Case I or Case III applies.
    In both cases, $\phi(w) = v_{t-1}$, so \ref{strong1} holds.
    \item[\scriptsize$\bullet$] If $x_1=w$ and $n'=n$, then again either Case I or Case III applies. Consider the vertex $u':=\phi^{-1}(v)$. Since $u'\in V(T'_{1})$, all of its neighbours are in $V(T'_1)\cup \{w\}$, and by our choice of $\phi_{0}$ and $\phi'$, $\phi(V(T'_1)\cup \{w\})\subseteq N_{G^2_T}(v_{t-1})\cup \{v_{t-1}\}$. We modify $\phi$ by letting $\phi(x_1):=v$ and $\phi(u'):=v_{t-1}$. Then $\phi$ is a $T$-admissible embedding of $T'$ into $G^2_T$ satisfying \ref{strong1}.
\end{enumerate}

This proves $A_{k,t}$. The lemma holds by induction.
\end{proof}

\subsection{Constructing $G_{T}^{2}$}

We now construct a $(4,1)$-tree $T$ on $n$ vertices such that $G^2_T$ contains at most $\frac{14}{5}n\log_2 n+\bigo(n)$ edges. This construction together with \Cref{strong} proves \Cref{thm2}. 
We construct this tree level by level, and assign a `type' to each vertex.

\begin{figure}[!h]
    \centering
    \begin{tikzpicture}[/tikz/xscale=1.3, /tikz/yscale=0.9]

    \begin{scope}[every node/.style={draw, shape = circle, fill = black, minimum size = 0.1cm, inner sep=1.7pt}]
        \node (01) at (0,3*1.3) {};
        \node (11) at (0,2*1.3) {};
        
        \node (21) at (-3*1.3,1*1.3) {};
        \node (22) at (-2*1.3,1*1.3) {};
        \node (23) at (-1*1.3,1*1.3) {};
        \node (24) at (0,1*1.3) {};
        \node (25) at (1*1.3,1*1.3) {};
        \node (26) at (2*1.3,1*1.3) {};
        \node (27) at (3*1.3,1*1.3) {};

        \node (31) at (-3*1.3,0) {};

        \node (32) at (-3*1.3+0.325*1,0) {};
        \node (33) at (-3*1.3+0.325*2,0) {};
        \node (34) at (-3*1.3+0.325*3,0) {};
        \node (35) at (-3*1.3+0.325*4,0) {};
        \node (36) at (-3*1.3+0.325*5,0) {};
        \node (37) at (-3*1.3+0.325*6,0) {};
        \node (38) at (-3*1.3+0.325*7,0) {};
        
        \node (39) at (-1*1.3,0) {};

        \node (310) at (-1*1.3+0.325*1,0) {};
        \node (311) at (-1*1.3+0.325*2,0) {};
        \node (312) at (-1*1.3+0.325*3,0) {};
        \node (313) at (-1*1.3+0.325*4,0) {};
        \node (314) at (-1*1.3+0.325*5,0) {};
        \node (315) at (-1*1.3+0.325*6,0) {};
        \node (316) at (-1*1.3+0.325*7,0) {};
        
        \node (317) at (1*1.3,0) {};

        \node (318) at (1*1.3+0.325*1,0) {};
        \node (319) at (1*1.3+0.325*2,0) {};
        \node (320) at (1*1.3+0.325*3,0) {};
        \node (321) at (1*1.3+0.325*4,0) {};
        \node (322) at (1*1.3+0.325*5,0) {};
        \node (323) at (1*1.3+0.325*6,0) {};
        \node (324) at (1*1.3+0.325*7,0) {};

        \node (325) at (3*1.3,0) {};
    \end{scope}

    \begin{scope}[every edge/.style={draw, thick}]
        \path [-] (01) edge (11);
        \path [-] (11) edge (21);
        \path [-] (11) edge (22);
        \path [-] (11) edge (23);
        \path [-] (11) edge (24);
        \path [-] (11) edge (25);
        \path [-] (11) edge (26);
        \path [-] (11) edge (27);

        \path [-] (21) edge (31);
        
        \path [-] (22) edge (32);
        \path [-] (22) edge (33);
        \path [-] (22) edge (34);
        \path [-] (22) edge (35);
        \path [-] (22) edge (36);
        \path [-] (22) edge (37);
        \path [-] (22) edge (38);

        \path [-] (23) edge (39);

        \path [-] (24) edge (310);
        \path [-] (24) edge (311);
        \path [-] (24) edge (312);
        \path [-] (24) edge (313);
        \path [-] (24) edge (314);
        \path [-] (24) edge (315);
        \path [-] (24) edge (316);

        \path [-] (25) edge (317);

        \path [-] (26) edge (318);
        \path [-] (26) edge (319);
        \path [-] (26) edge (320);
        \path [-] (26) edge (321);
        \path [-] (26) edge (322);
        \path [-] (26) edge (323);
        \path [-] (26) edge (324);

        \path [-] (27) edge (325);
    \end{scope}

    \begin{scope}
        \node (01) at (0+0.35,3*1.3) {$1$};
        \node (11) at (0+0.35,2*1.3+0.2) {$2$};
        
        \node (21) at (-3*1.3-0.35,1*1.3) {$1$};
        \node (22) at (-2*1.3-0.35,1*1.3) {$2$};
        \node (23) at (-1*1.3-0.35,1*1.3) {$1$};
        \node (24) at (0+0.35,1*1.3) {$2$};
        \node (25) at (1*1.3+0.35,1*1.3) {$1$};
        \node (26) at (2*1.3+0.35,1*1.3) {$2$};
        \node (27) at (3*1.3+0.35,1*1.3) {$1$};

        \node (31) at (-3*1.3,0-0.4) {$2$};

        \node (32) at (-3*1.3+0.325*1,0-0.4) {$1$};
        \node (33) at (-3*1.3+0.325*2,0-0.4) {$2$};
        \node (34) at (-3*1.3+0.325*3,0-0.4) {$1$};
        \node (35) at (-3*1.3+0.325*4,0-0.4) {$2$};
        \node (36) at (-3*1.3+0.325*5,0-0.4) {$1$};
        \node (37) at (-3*1.3+0.325*6,0-0.4) {$2$};
        \node (38) at (-3*1.3+0.325*7,0-0.4) {$1$};
        
        \node (39) at (-1*1.3,0-0.4) {$2$};

        \node (310) at (-1*1.3+0.325*1,0-0.4) {$1$};
        \node (311) at (-1*1.3+0.325*2,0-0.4) {$2$};
        \node (312) at (-1*1.3+0.325*3,0-0.4) {$1$};
        \node (313) at (-1*1.3+0.325*4,0-0.4) {$2$};
        \node (314) at (-1*1.3+0.325*5,0-0.4) {$1$};
        \node (315) at (-1*1.3+0.325*6,0-0.4) {$2$};
        \node (316) at (-1*1.3+0.325*7,0-0.4) {$1$};
        
        \node (317) at (1*1.3,0-0.4) {$2$};

        \node (318) at (1*1.3+0.325*1,0-0.4) {$1$};
        \node (319) at (1*1.3+0.325*2,0-0.4) {$2$};
        \node (320) at (1*1.3+0.325*3,0-0.4) {$1$};
        \node (321) at (1*1.3+0.325*4,0-0.4) {$2$};
        \node (322) at (1*1.3+0.325*5,0-0.4) {$1$};
        \node (323) at (1*1.3+0.325*6,0-0.4) {$2$};
        \node (324) at (1*1.3+0.325*7,0-0.4) {$1$};

        \node (325) at (3*1.3,0-0.4) {$2$};
    \end{scope}
    \end{tikzpicture}
    
    \caption{$T_{k}$ when $k=3$.}
\end{figure}
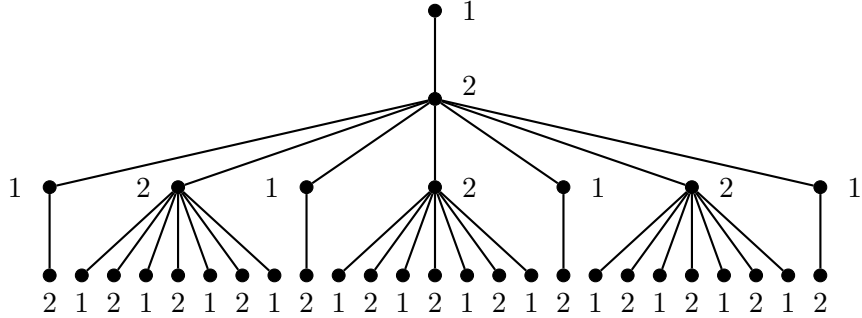



For a given integer $k\geq0$, we define a tree $T_{k}$ of level $k$ as follows. Note that the children of each vertex will be ordered from left to right.
Start with a root vertex $v$ of type $1$. For every vertex $x$ on level $\ell<k$ of type $p\in\{1, 2\}$, we add its children as follows.
\begin{enumerate}
    \item[\scriptsize$\bullet$] If $p=1$, then add exactly one child of type $2$.
    \item[\scriptsize$\bullet$] If $p=2$, then add exactly seven children of type $1, 2, 1, 2, 1, 2, 1$ from left to right.
\end{enumerate}
Repeating this for all vertices at levels less than $k$ yields our desired tree $T_{k}.$
For each $\ell\leq k$, let $x^{\ell}_1,\dots, x^{\ell}_{t_{\ell}}$ be the vertices on level $\ell$ ordered from left to right. 

From this definition, it is evident that all vertices of the same type on the same level have the same number of descendants. 
We define the sequence $a_{\ell,p}=a^{k}_{\ell,p}$ to be $\nu_{T_{k}}(x)$, the number of descendants of a vertex $x$ (the vertex itself is also counted) at level $\ell$ of type $p$. We will omit the superscript $k$ as it is clear from the context. Then the above definition yields the following recursive relation on the sequence $a_{\ell,p}$.
\begin{align}\label{eq: definition of a}
a_{\ell, p}=
\begin{cases}
1 & \text{ if } \ell=k\;, \\
a_{\ell+1, 2}+1 & \text{ if } \ell <k\text{ and }p=1\;,\\
4a_{\ell+1, 1}+3a_{\ell+1, 2}+1 & \text{ if } \ell <k\text{ and }p=2\;.\end{cases}
\end{align}
Note that there are no vertices of type $2$ at level $0$ and no vertices of type $1$ at level $1$, but the numbers $a_{1,1}$ and $a_{0,2}$ can still be defined by the recurrence relation described above. From this relation, one can deduce the following properties:
\begin{enumerate}[label=\textup{\textbf{R\arabic{enumi}}},ref=\textup{\textbf{R\arabic{enumi}}}]
    \item\label{R1} $a_{\ell, 2}-1=4(a_{\ell, 1}-1)$ if $k-\ell$ is even and  $a_{\ell, 2}=4a_{\ell, 1}$ if $k-\ell$ is odd. In particular, $4a_{\ell, 1}-3\leq a_{\ell, 2}\leq4a_{\ell, 1}$.
    
    \item\label{R2} For $\ell\in\{0,\ldots,k\}$, the types for the vertices $x^{\ell}_1,\dots, x^{\ell}_{t_{\ell}}$ in level $\ell$ (from left to right) are $1,2,\dots,2,1$ when $\ell$ is even, and $2,1,\dots,1,2$ when $\ell$ is odd. In particular, $t_{\ell}$ is always odd.

    \item \label{R3} $|T_{k}|=4|T_{k-2}|+3|T_{k-1}|-1$, which implies $\frac{1}{2}+2^{2k-1}\leq|T_{k}|\leq2^{2k+1}-1$ for all $k\geq 2$. In particular, $k\leq \frac{1}{2}+ \frac{1}{2}\log_2|T_{k}|$ holds for all $k\geq 0$.
    
    \item \label{R4} For each $\ell\in [k]$, $\sum_{\ell\leq \ell'\leq k} a_{\ell',1} \leq 2 a_{\ell,1}$.
\end{enumerate}
Indeed, \ref{R1} and \ref{R2} are routine to check. \ref{R3} can be seen to hold by the recursive structure of $T_{k}$. \ref{R4} is implied by the relations:
$a_{\ell,1} = a_{\ell+1,2}+1 \geq 4 a_{\ell+1,1}-2 \geq 2a_{\ell+1,1}$ for all $\ell<k$, where the first relation holds by \eqref{eq: definition of a}, and the second relation holds by \ref{R1}. 


One can verify that $T_k$ is a $(4,1)$-tree: \ref{R2} can be used to verify  \ref{T1} from \Cref{def: Ks tree}, \ref{R1} can be used to verify \ref{T2}. Furthermore, \ref{T3} and \ref{T4} can be seen to hold by the structure of $T_{k}$.



We now estimate the number of edges in $G_{T_k}^{2}$. We begin by first estimating the number of arcs in the digraph $\overrightarrow{G}_{T_{k}}$ which will provide an upper bound on $e(G_{T_k})$.



Let $x\in V(T_{k})$ be a vertex of level $\ell\geq 2$. If $x$ is of type $1$, then the child of $x$ sends at most $\nu(l(x))+\nu(x)$ outgoing arcs. If $x$ is of type $2$, then each of the seven children of $x$ (from left to right) sends outgoing arcs to its own descendants and the descendants of its left-siblings. Thus, each of the seven children of $x$ (from left to right) send at most
$$a_{\ell+1,1},\;\; a_{\ell+1,1}+a_{\ell+1,2},
\;\; 2 a_{\ell+1,1}+a_{\ell+1,2},\;\; 2a_{\ell+1,1}+2a_{\ell+1,2},
\;\; \dots,\;\; 4a_{\ell+1,1}+3a_{\ell+1,2}$$ 

outgoing arcs into $D_{T_{k}}[x]$, which totals to $16 a_{\ell+1,1}+ 12 a_{\ell+1,2}\leq 4\nu(x)$ arcs.
Additionally, each of them sends $\nu(l(x))$ outgoing arcs into $D_T[l(x)]$. Therefore, in total, the children of $x$ send at most $7\nu(l(x))+4\nu(x)$ arcs.
Hence, the number of arcs sent from the vertices of level $\ell+1$ is at most (where $\nu(x^{\ell+1}_0):=0$)
\begin{align*}
    \sum_{x^{\ell+1}_i  \text{ is type } 2} \left(7\nu(x^{\ell+1}_{i-1})+ 4\nu(x^{\ell+1}_i)\right) +& \sum_{x^{\ell+1}_i \text{ is type } 1} \left(\nu(x^{\ell+1}_{i-1})+\nu(x^{\ell+1}_i)\right)\\
    &\quad<\sum_{x^{\ell+1}_i \text{ is type }2} 5 \nu(x^{\ell+1}_{i}) + \sum_{x^{\ell+1}_i  \text{ is type } 1} 8\nu(x^{\ell+1}_{i})
\end{align*}
because vertices of type $1$ and $2$ alternate by \ref{R2}. Then,
\begin{align*}
    \sum_{x^{\ell+1}_i \text{ is type } 2} 5 \nu(x^{\ell+1}_{i}) + \sum_{x^{\ell+1}_i \text{ is type } 1} 8\nu(x^{\ell+1}_{i})&=5\cdot\sum_{i\in[t_{\ell+1}]}\nu(x^{\ell+1}_{i})+3\cdot\sum_{x^{\ell+1}_i \text{ is type } 1} a_{\ell+1,1}\\
    &<\frac{28}{5}\sum_{i\in[t_{\ell+1}]}\nu(x^{\ell+1}_{i})+3t_{\ell+1}+3a_{\ell+1, 1}\\
    &\leq\frac{28}{5}|T_{k}|+3t_{\ell+1}+3a_{\ell+1, 1}\;,
\end{align*} 
where the second inequality follows from \ref{R1}, 
and the last inequality holds since $\sum_{i\in [t_{\ell+1}]}\nu(x^{\ell+1}_i)\leq|T_{k}|-1$ when $\ell\geq 2$. Since $T_{k}$ has $k$ levels, and the nine vertices of level at most $2$ each send at most $|T_{k}|$ outgoing arcs, 
\begin{align}\label{eq: edges in Trk}
    e(\overrightarrow{G}_{T_{k}}) \leq 9|T_{k}| + \frac{28}{5}k|T_{k}|+3\sum_{\ell=2}^{k-1}t_{\ell+1}+3\sum_{\ell=2}^{k-1}a_{\ell+1, 1}< \frac{14}{5} |T_{k}| \log_2 |T_{k}| + 18 |T_{k}|\;,
\end{align} 
where the last inequality holds by \ref{R3}, $\sum_{\ell=2}^{k-1}a_{\ell+1, 1}\leq 2a_{3,1}<|T_k|$ by \ref{R4}, and because $\sum_{\ell=2}^{k-1}t_{\ell+1}\leq |T_k|$. 

Now, for a given sufficiently large $n$, let $k$ be the integer satisfying $|T_{k-1}|<n\leq |T_{k}|$. Note, by \ref{R3} it follows that $|T_{k}|\leq16n$. Consider the $n$-vertex $T_{k}$-admissible set $U$. Let $G$ be the subgraph of $G_{T_{k}}$ induced by $U$. We estimate the number of edges in $G$.

Let $X\subseteq U$ be the set of vertices $x$ such that $D_{T_{k}}[x] \subseteq U$ while $D_{T_{k}}[x^*]\not\subseteq U$ for its parent $x^*$.  Define $\define{$D[X]$}:=\bigcup_{x\in X} D_{T_{k}}[x]$. The following claim is useful to bound the size of $U\setminus D[X]$.
\begin{lemma}\label{lem: U-D[X]}
The set $U\setminus D[X]$ contains at most one vertex from each level. 
\end{lemma}
\begin{proof}
Consider the DFS ordering $x_1,\dots, x_{|T_{k}|}$ of $T_{k}$. Since $U$ is $T_{k}$-admissible, $U=\{x_1,\dots, x_{n}\}$. For a contradiction, assume that there exists $x_{i}, x_{j} \not\in D[X]$ with $i<j\leq n$ and $L(x_{i})=L(x_{j})$. Observe that the DFS preorder traversal visits $x_j$ only after visiting all of the descendants of $x_i$. Therefore, $D_{T_{k}}[x_i]\subseteq U$, which implies $x_i \in D[X]$; a contradiction. This proves the lemma.
\end{proof}

In order to estimate $e(G)$, we begin by counting the arcs in $\overrightarrow{G}_{T_{k}}[U]$. 
First, consider the outgoing arcs from $U\setminus D[X]$. For each vertex $u\in U\setminus D[X]$ of level $\ell\geq 2$, $u$ sends at most $\nu(u^*)+ \nu(l(u^*))$ arcs into $U$. By \ref{R2}, $\nu(u^*)+ \nu(l(u^*)) \leq a_{\ell-1,1}+ a_{\ell-1,2}$. By \Cref{lem: U-D[X]}, $U\setminus D[X]$ contains at most one vertex from each level $\ell$, so the number of outgoing arcs from $U\setminus D[X]$ is at most 
\begin{align}\label{eq: edges in U1}
    2n+\sum_{2\leq \ell\leq k} \left(a_{\ell-1,1}+a_{\ell-1,2} \right)\leq 2n + 5 \sum_{2\leq \ell \leq k} a_{\ell-1,1}\leq2n + 10|T_{k}|\leq 162n\;,
\end{align}
where the $2n$ term counts the number of arcs from the two vertices at level at most one, the first inequality follows from \ref{R1}, and the second inequality follows from \ref{R4} since $a_{1,1}\leq a_{1,2}\leq|T_{k}|$.


Second, for each $x\in X$, we count the arcs with both ends in $D[x]=D_{T_{k}}[x]$. Let $\ell$ be the level of $x$.
If $x$ is type $1$, then the subgraph induced by $D[x]$ is isomorphic to $T_{k-\ell}$. If $x$ is type $2$, then the subgraph induced by $D[x]$ is isomorphic to a copy of $T_{k-\ell+1}$ minus its root. Thus, \eqref{eq: edges in Trk} implies that the number of arcs with both ends in $D[x]$ is at most 
$\frac{14}{5}(|D[x]|+1)\log_2 (|D[x]|+1) + 18(|D[x]|+1) \leq \frac{14}{5}|D[x]|\log_2 |D[x]| + 50|D[x]|.$

Third, for each $x\in X$, we count the arcs from $D[x]$ to $U\setminus D[x]$. Let $\ell$ be the level of $x$. Let $w\in D[x]$. If $w=x$, then since $w$ has at most $6$ left-siblings, by \ref{G2} and \ref{G3}, $w$ sends at most $a_{\ell-1, 2}+6a_{\ell, 2}$ arcs. Thus by \ref{R1}, $w$ sends at most $a_{\ell-1,2}+6a_{\ell, 2} \leq 4a_{\ell-1, 1}+6a_{\ell, 2}\leq 4( a_{\ell,2}+1)+6a_{\ell, 2} \leq 
44a_{\ell, 1}\leq 44|D[x]|$ arcs to $U\setminus D[x]$. If $w\neq x$, then the arcs from $w$ to $U\setminus D[x]$ must be arcs from $w$ to $D[l(w^*)]$ where $l(w^*)\notin D[x]$. This occurs only when $w^*$ is the left-most vertex in $D[x]$ of a level $\ell' \geq \ell$. Since $w^{*}$ has at most $7$ children, $w$ sends at most $7a_{\ell', 2}$ arcs. Thus, the total number of arcs from $D[x]$ to $U\setminus D[x]$ is at most
\begin{align}\label{eq: edges in U2}
    44|D[x]| + \sum_{\ell'\geq \ell} 7 a_{\ell', 2}\leq 44|D[x]|+ 28 \cdot \sum_{\ell'\geq \ell} a_{\ell', 1}\leq 44|D[x]|+56 a_{\ell, 1}\leq 100 |D[x]|\;,
\end{align}
where the first inequality holds by \ref{R1}, the second inequality holds by \ref{R4}, and the third inequality holds since $a_{\ell, 1}\leq|D[x]|$. Hence, the total number of arcs in $\overrightarrow{G}_{T_{k}}[U]$ is at most
\begin{align*}
    \#(\text{arcs from }U\setminus D[X]\text{ to }G)&+\sum_{x\in X}\Big[\#(\text{arcs from }D[x]\text{ to }D[x])+\#(\text{arcs from }D[x]\text{ to }U\setminus D[x])\Big]\\
    &\leq 162n + \sum_{x\in X} \left(\frac{14}{5}|D[x]| \log_2|D[x]| + 50 |D[x]|\right) + \sum_{x\in X} 100 |D[x]| \\
    &<  \frac{14}{5} n \log_2 n + 400 n\;,
\end{align*}
where the first inequality follows from \labelcref{eq: edges in U1,eq: edges in U2}, and the second inequality holds since $|\sum_{x\in X}|D[x]||\leq n$. Therefore, $e(G)\leq \frac{14}{5} n \log_2 n + 400 n$.

By \Cref{strong}, $G^2_{T_k}[U]$ is an $n$-vertex $\mathcal{T}_{n}$-universal graph. Furthermore, $G^2_{T_k}[U]$ has at most $3\cdot 7^2 n \leq 200 n$ more edges than $G$ from \Cref{newdef}~\ref{G4} (this is because every vertex in $T_k$ has at most $7$ children). Hence, $e(G^2_{T_k}[U])\leq e(G)+200 n\leq\frac{14}{5} n \log_2 n + 600 n$. We conclude that $G^2_{T_{k}}[U]$ is the desired universal graph, which proves \Cref{thm2}.

\section{Interval-Universality}\label{sec:6}

For any $(4,1)$-tree $T$, \Cref{strong} yields an embedding $\phi$ of a given tree $T'$ with $|T'|\leq |T|$ into the universal graph $G^2_T$. Moreover, it ensures that $V(G^2_T)\setminus\phi(V(T'))$ is a $T$-admissible set. Therefore, it follows that for any interval $U=\{x_i,x_{i+1},\dots, x_j\}$ in the DFS preorder traversal $x_1,\dots, x_{|T|}$, the induced subgraph $G^2_T[U]$ is $\mathcal{T}_{|U|}$-universal. This motivates the following new definition.


\begin{definition}
An $n$-vertex graph $G$ is \define{interval-universal} (for trees) if there exists an ordering $x_1,\dots, x_n$ of the vertices of $G$ such that for all integers $i,m\geq0$, the induced subgraph $G[\{x_{i+1},\dots, x_{i+m}\}]$ is $\mathcal{T}_{m}$-universal.
\end{definition}

Define \define{$s^{\rm int}(n)$} to be the minimum number of edges in an $n$-vertex graph that is interval-universal for trees. Then, 
\Cref{strong} implies that for any $(4,1)$-tree $T$, the $2$-generated graph $G^{2}_T$ is interval-universal with respect to the DFS preorder traversal. Hence, we have the following bounds relating $s^{\rm int}(n)$, $s^{*}(n)$ and $s(n)$.
\begin{align}\label{intervaluniversality: ineq 1}
n\ln{n}-\bigo(n)\leq s^{*}(n)\leq s(n)\leq s^{\rm int}(n)\leq\frac{14}{5}n\log_2{n}+\bigo(n)\;.
\end{align}
Improving the lower bounds for $s^{*}(n)$ and $s(n)$ seems very difficult. However, for $s^{\rm int}(n)$, we can improve the lower bound.

\begin{theorem}\label{thm3}
$s^{\rm int}(n)\geq\left(1-o(1)\right)n\log_2{n}$.
\end{theorem}
\begin{proof}
We proceed by induction on $n$ with the hypothesis that $s^{\rm int}(n) \geq n \log_2 n - 4n\sqrt{\log_2 n}$. When $n=1$, $n\log_2 n-4n\sqrt{\log_2 n}=0$, so the result holds trivially. Now assume that $n > 1$ and $s^{\rm int}(n') \geq n' \log_2 n' - 4n'\sqrt{\log_2 n'}$ for all $n'<n$. Let $\psi(x)=x \log_2 x - 4x\sqrt{\log_2 x}$. Since the second derivative of $\psi(x)$ is positive when $x>1$, $\psi(x)$ is convex.

Let $G$ be an interval-universal graph for trees, and let $x_1,\dots, x_n$ be the associated ordering of $V(G)$. 
Since the $n$-vertex star is contained in $G$, there is a vertex $x_t$ of degree $n-1$ in $G$.
Since $G[ \{x_1,\dots, x_{t-1}\}]$ and $G[\{x_{t+1},\dots, x_n\}]$ are both interval-universal, by induction,
\begin{align*}
    e(G)\geq \psi(t-1)+\psi(n-t)+n-1&\geq 2\psi\Big(\frac{n-1}{2}\Big)+n-1\\
    &=  (n-1)\log_2\Big(\frac{n-1}{2}\Big) - 4(n-1)\sqrt{\log_2\Big(\frac{n-1}{2}\Big)} + n-1\\
    &\geq n\log_2 n  - 4n \sqrt{\log_2 n}\;,
\end{align*}
where the second inequality holds by Jensen's inequality, and the last inequality holds since $(n-1)\log_2(\frac{n-1}{2}) \geq n\log_2 n - (n-1) -\log_2 n -2$, $4(n-1)\sqrt{\log_2(\frac{n-1}{2})} \leq 4n\sqrt{\log_2 n} - \frac{2n}{\sqrt{\log_2{n}}}$ and $\frac{2n}{\sqrt{\log_2 n}} \geq \log_2 n +2$. Hence, $s^{\rm int}(n)\geq (1-o(1)) n \log_2 n$.
\end{proof}

Note that the gap between the lower bound in \Cref{thm3} and the upper bound in \labelcref{intervaluniversality: ineq 1} is at most the multiplicative constant $\frac{14}{5}$. A natural problem is to improve the upper bound.
One idea is to increase the imbalance $K$ to a value larger than $4$ as explained in \Cref{sec: overview} without adding edges from $u$ to the descendants of its nearest-right cousin. This will reduce the number of edges, however, proving that such a graph is universal will need new ideas.

Note that the strategy of generating a graph $G^2_T$ from a given tree $T$ and using admissibility in an induction argument yields a bound on $s^{\rm int}(n)$. Thus, just as \Cref{thm3} distinguishes $s^{\rm int}(n)$ from $s^*(n)$ or $s(n)$, some new ideas that only apply to $s^*(n)$ or $s(n)$, but not necessarily to $s^{\rm int}(n)$, might be helpful to further improve upper bounds on $s^*(n)$ and $s(n)$. 

\begin{problem}
Is there a constant $c<\frac{14}{5\ln{2}}$ such that $s^*(n)\leq (c+o(1))n\ln{n}$? 
\end{problem}

\subsection*{Acknowledgement}

Thanks to Carla Groenland, who identified some technical issues in an earlier version of the paper, which are now fixed. 

{
\fontsize{11pt}{12pt}
\selectfont
\bibliographystyle{plainnat}
\bibliography{references} 

\begin{thebibliography}{23}
\providecommand{\natexlab}[1]{#1}
\providecommand{\url}[1]{\texttt{#1}}
\expandafter\ifx\csname urlstyle\endcsname\relax
  \providecommand{\doi}[1]{doi: #1}\else
  \providecommand{\doi}{doi: \begingroup \urlstyle{rm}\Url}\fi

\bibitem[Alon and Asodi(2002)]{alon}
N.~Alon and V.~Asodi.
\newblock Sparse universal graphs.
\newblock \emph{Journal of Computational and Applied Mathematics}, 142\penalty0 (1):\penalty0 1--11, 2002.

\bibitem[Alon and Capalbo(2007)]{alonsecond}
N.~Alon and M.~Capalbo.
\newblock Sparse universal graphs for bounded-degree graphs.
\newblock \emph{Random Structures \& Algorithms}, 31\penalty0 (2):\penalty0 123--133, 2007.

\bibitem[Alon et~al.(2001)Alon, Capalbo, Kohayakawa, Rödl, Ruciński, and Szemerédi]{szemeredi}
N.~Alon, M.~Capalbo, Y.~Kohayakawa, V.~Rödl, A.~Ruciński, and E.~Szemerédi.
\newblock Near-optimum universal graphs for graphs with bounded degrees.
\newblock \emph{Approximation, Randomization, and Combinatorial Optimization: Algorithms and Techniques}, 5:\penalty0 170--180, 2001.

\bibitem[Alon et~al.(2007)Alon, Krivelevich, and Sudakov]{alonthird}
N.~Alon, M.~Krivelevich, and B.~Sudakov.
\newblock Embedding nearly-spanning bounded degree trees.
\newblock \emph{Combinatorica}, 27\penalty0 (6):\penalty0 170--180, 2007.

\bibitem[Babai et~al.(1982)Babai, Chung, Erd{\"o}s, Graham, and Spencer]{babai}
L.~Babai, F.~Chung, P.~Erd{\"o}s, R.~Graham, and J.~Spencer.
\newblock On graphs which contain all sparse graphs.
\newblock \emph{North-Holland Mathematics Studies}, 60:\penalty0 21--26, 1982.

\bibitem[Bergold et~al.(2024)Bergold, Iršič, Lauff, Orthaber, Scheucher, and Wesolek]{bergold2024}
H.~Bergold, V.~Iršič, R.~Lauff, J.~Orthaber, M.~Scheucher, and A.~Wesolek.
\newblock Subgraph-universal planar graphs for trees.
\newblock \emph{arXiv:2409.01678}, 2024.
\newblock URL \url{https://arxiv.org/abs/2409.01678}.

\bibitem[Bhatt et~al.(1989)Bhatt, Chung, Leighton, and Rosenberg]{bhatt}
S.~Bhatt, F.~Chung, F~Leighton, and A.~Rosenberg.
\newblock Universal graphs for bounded-degree trees and planar graphs.
\newblock \emph{SIAM Journal on Discrete Mathematics}, 2\penalty0 (2):\penalty0 145--155, 1989.

\bibitem[Bucić et~al.(2021)Bucić, Draganić, and Sudakov]{bucic}
M.~Bucić, N.~Draganić, and B.~Sudakov.
\newblock Universal and unavoidable graphs.
\newblock \emph{Combinatorics, Probability and Computing}, 30\penalty0 (6):\penalty0 942–955, 2021.

\bibitem[Capalbo(2002)]{capalbo}
M.~Capalbo.
\newblock Small universal graphs for bounded-degree planar graphs.
\newblock \emph{Combinatorica}, 22\penalty0 (3):\penalty0 345--359, 2002.

\bibitem[Capalbo and Kosaraju(1999)]{capalbosecond}
M.~Capalbo and S.~Kosaraju.
\newblock Small universal graphs.
\newblock \emph{Proceedings of the Thirty-First Annual ACM Symposium on Theory of Computing}, page 741–749, 1999.

\bibitem[Chung(1990)]{chung}
F.~Chung.
\newblock Separator theorems and their applications.
\newblock In \emph{Paths, flows, and VLSI-layout}, volume~9 of \emph{Algorithms Comb.}, pages 17--34. Springer, 1990.

\bibitem[Chung and Graham(1978)]{grahamfirst}
F.~Chung and R.~Graham.
\newblock On graphs which contain all small trees.
\newblock \emph{Journal of Combinatorial Theory, Series B}, 24:\penalty0 14--23, 1978.

\bibitem[Chung and Graham(1983)]{grahamthird}
F.~Chung and R.~Graham.
\newblock On universal graphs for spanning trees.
\newblock \emph{Journal of the London Mathematical Society}, s2-27\penalty0 (2):\penalty0 203--211, 1983.

\bibitem[Chung and Graham(2006)]{grahamsecond}
F.~Chung and R.~Graham.
\newblock On universal graphs.
\newblock \emph{Annals of the New York Academy of Sciences}, 319\penalty0 (1):\penalty0 136--140, 2006.

\bibitem[Chung et~al.(1978)Chung, Graham, and Pippenger]{pipenger}
F.~Chung, R.~Graham, and N.~Pippenger.
\newblock On graphs which contain all small trees, \uppercase\expandafter{\romannumeral2}.
\newblock \emph{Colloquia Mathematica Societatis János Bolyai}, 18:\penalty0 213--223, 1978.

\bibitem[Esperet et~al.(2023)Esperet, Joret, and Morin]{joret}
L.~Esperet, G.~Joret, and P.~Morin.
\newblock Sparse universal graphs for planarity.
\newblock \emph{Journal of the London Mathematical Society}, 108\penalty0 (4):\penalty0 1333--1357, 2023.

\bibitem[Frati et~al.(2023)Frati, Hoffmann, and Tóth]{frati}
F.~Frati, M.~Hoffmann, and C.~Tóth.
\newblock Universal geometric graphs.
\newblock \emph{Combinatorics, Probability and Computing}, 32\penalty0 (5):\penalty0 742–761, 2023.

\bibitem[Gy\H{o}ri et~al.(2025)Gy\H{o}ri, Li, Salia, and Tompkins]{gyori}
E.~Gy\H{o}ri, B.~Li, N.~Salia, and C.~Tompkins.
\newblock A note on universal graphs for spanning trees.
\newblock \emph{Discrete Applied Mathematics}, 362:\penalty0 146–147, 2025.

\bibitem[Joret(2020)]{JoretTalk}
G.~Joret.
\newblock Sparse universal graphs for planarity.
\newblock \emph{Oxford Discrete Maths and Probability Seminar}, 2020.
\newblock URL \url{https://www.youtube.com/watch?v=6-vz5GmBrSA}.

\bibitem[Kaul and Wood(2025)]{Kaul}
N.~Kaul and D.~R. Wood.
\newblock On universal graphs for trees and treewidth $k$ graphs.
\newblock \emph{arXiv preprint arXiv:2508.0335}, 2025.

\bibitem[Kim and Kim(2025)]{KimKim}
J.~Kim and M.~Kim.
\newblock On the size of universal graphs for spanning trees.
\newblock \emph{arXiv preprint arXiv:2508.19032}, 2025.

\bibitem[Montgomery(2019)]{montgomery}
R.~Montgomery.
\newblock Spanning trees in random graphs.
\newblock \emph{Advances in Mathematics}, 356:\penalty0 \#106793, 2019.

\bibitem[Wood(2025)]{wood}
D.~R. Wood.
\newblock Tree decompositions with small width, spread, order and degree.
\newblock \emph{arXiv preprint arXiv:2509.01140}, 2025.

\end{thebibliography}
}

\end{document}